\documentclass[fleqn,final,3p,11pt]{elsarticle}

\usepackage{bbm}
\usepackage{lineno,hyperref}
\hypersetup{pdfauthor=author}
\usepackage{amsmath,amsthm,amssymb}
\usepackage{dirtytalk}
\usepackage{comment}

\usepackage{siunitx}
\usepackage{xcolor}
\usepackage{booktabs,colortbl, array}

\definecolor{rulecolor}{RGB}{0,71,171}
\definecolor{tableheadcolor}{gray}{0.92}
\usepackage{rotating}

\usepackage{graphicx}
\usepackage{multirow}
\usepackage{colortbl}%
\usepackage{booktabs}
\definecolor{Gray}{gray}{0.925}
\usepackage{enumitem,booktabs,cfr-lm}
\usepackage[referable]{threeparttablex}
\renewlist{tablenotes}{enumerate}{1}
\makeatletter
\setlist[tablenotes]{label=\tnote{\alph*},ref=\alph*,itemsep=\z@,topsep=\z@skip,partopsep=\z@skip,parsep=\z@,itemindent=\z@,labelindent=\tabcolsep,labelsep=.2em,leftmargin=*,align=left,before={\footnotesize}}
\makeatother
\usepackage{placeins,diagbox}
\usepackage{booktabs} % For pretty tables
\usepackage{caption} % For caption spacing
\usepackage{subcaption} % For sub-figures

\usepackage{graphicx}
\usepackage[all]{nowidow}
\usepackage[utf8]{inputenc}
\usepackage{tikz}
\usetikzlibrary{er,positioning,bayesnet}
\usepackage{multicol}
\usepackage{algpseudocode,algorithm,algorithmicx}

 \newtheorem{thm}{Theorem}[section]

 \numberwithin{equation}{section}
\modulolinenumbers[5]
\begin{document}

\begin{frontmatter}

%% Title, authors and addresses

%% use the tnoteref command within \title for footnotes;
%% use the tnotetext command for theassociated footnote;
%% use the fnref command within \author or \address for footnotes;
%% use the fntext command for theassociated footnote;
%% use the corref command within \author for corresponding author footnotes;
%% use the cortext command for theassociated footnote;
%% use the ead command for the email address,
%% and the form \ead[url] for the home page:
%% \title{Title\tnoteref{label1}}
%% \tnotetext[label1]{}
%% \author{Name\corref{cor1}\fnref{label2}}
%% \ead{email address}
%% \ead[url]{home page}
%% \fntext[label2]{}
%% \cortext[cor1]{}
%% \affiliation{organization={},
%%             addressline={},
%%             city={},
%%             postcode={},
%%             state={},
%%             country={}}
%% \fntext[label3]{}

\title{On a Stable Method for Option Pricing: Discontinuous Petrov-Galerkin Method for Option Pricing and Sensitivity Analysis}

%% use optional labels to link authors explicitly to addresses:
%% \author[label1,label2]{}
%% \affiliation[label1]{organization={},
%%             addressline={},
%%             city={},
%%             postcode={},
%%             state={},
%%             country={}}
%%
%% \affiliation[label2]{organization={},
%%             addressline={},
%%             city={},
%%             postcode={},
%%             state={},
%%             country={}}

\author[inst1]{Davood Damircheli}

\address[label1]{Department and Organization, Mississippi State University,\\ P.O. Box 39759, Starkville, USA}

\begin{abstract}
The discontinuous Petrov–Galerkin (DPG) methodology of Demkowicz and Gopalakrishnan introduced in \cite{demkowicz2010class} has been widely used for problems in computational mechanics. In this investigation, we propose the DPG method for option pricing and sensitivity analysis under the basic Black-Scholes model. In this investigation, primal and ultraweak formulation of the DPG method is derived for Vanilla options, American options, Asian options, and Barrier options. A wide range of standard numerical experiments is conducted to examine the convergence, stability, and efficiency of the proposed method for each one of the options separately. Besides, a C++ high-performance (HPC) code for option pricing with the DPG method is developed which is available to the public to customize it for option pricing problems or other related problems.
\end{abstract}
%%Graphical abstract
%\begin{graphicalabstract}
%\includegraphics{grabs}
%\end{graphicalabstract}

%%Research highlights
%\begin{highlights}
%\item Research highlight 1
%\item Research highlight 2
%\end{highlights}

\begin{keyword}
%% keywords here, in the form: keyword \sep keyword
Discontinuous Petrov–Galerkin \sep Quantitative Finance \sep Vanilla and Exotic Options\sep High performance Programming 
%% PACS codes here, in the form: \PACS code \sep code
%\PACS 0000 \sep 1111
%% MSC codes here, in the form: \MSC code \sep code
%% or \MSC[2008] code \sep code (2000 is the default)
%\MSC 0000 \sep 1111
\end{keyword}

\end{frontmatter}

%% \linenumbers

%% main text
\section{Introduction}
Since 1970, Black and Scholes \cite{black1973pricing} and Merton \cite{merton1973theory} discover the pioneering option pricing formula and received the Nobel Prize in Economics \cite{ferreyra1998mathematics}, known as the classical Black–Scholes (or BS–Merton) options pricing model, this formula has widespread been attentive for academia and practitioners. Options as an important element of the financial derivative market are hedging tools for investors to devise risk-protected strategies against fluctuations in the price of the underlying assets. Additionally, Greeks, the sensitivity of the option price with respect to the different parameters, are another standard gauge for market makers to design the best hedge for their position.   

However, the analytical solution for pricing these financial instruments except in some special cases is not obtainable. Among those instruments, Exotic options which are path-dependent financial derivatives, are very challenging problems to value. American-style options, Asian-type options, and Barrier options are the most appealing example of these families of hedging devices for traders, where they are  too complex problems to price analytically.

Admittedly, developing efficient numerical methods for option pricing as a natural remedy began shortly after developing the B-S model in 1970. To name some of the most widely used methods, one can cite analytical approximation \cite{barone1987efficient,geske1984american}, stochastic mesh method \cite{broadie2004stochastic}, Monte Carlo method \cite{boyle1997monte, boyle1977options, acworth1998comparison}, the lattice-based method for the corresponding partial differential equations like finite element, and finite difference methods \cite{achdou2005computational, chiarella2014numerical, duffy2013finite, seydel2006tools, tavella2000pricing}, and mesh-free methods \cite{kim2014option, fasshauer2004using,bastani2013radial} can be mentioned.

Among the aforementioned numerical methods, the weighted residual methods (or Galerkin methods) method has always been the center of interest from community \cite{seydel2006tools, achdou2005computational} due to the undoubted merit in the numerical solution of differential equations. Galerkin's methods benefit from the variety of advantages, having an elaborate and apprehended theory on prior and posterior error estimates, and significant flexibility for non-rectangle domains to just name a few. Thus, in the quantitative finance field, path-dependent options can use the well-understood error estimates to adapt a refiner mesh in their domain where needed. American options near the optimal exercise boundary or multi-factor options \cite{seydel2006tools, achdou2005computational} with the complicated domains can be a good example where the strength of variational methods can be simply exploited.

The state-of-art Black-Scholes PDE is a time-dependent parabolic partial differential equation that can be classified as a convection-diffusion equation. It is well-known in the literature (\cite{ern2004theory}, \cite{strang1974analysis}, and \cite{douglas1982numerical} references therein) that this family of problems can be numerically unstable once the coercivity violates due largely to the small coefficient of second order differential operator. This instability can emerge as a loss of accuracy or oscillatory behavior of the solution.

Since the discontinuous Petrov-Galerkin with optimal test space (DPG) method developed by Demkowicz and Golapalakrishnan \cite{demkowicz2010class}, it has been widely used for the numerical solution of Differential equation (\cite{mustapha2014discontinuous,roberts2015discontinuous,ernesti2019space,ellis2016space,fuhrer2017time,roberts2021time}, and the references therein) consist of convection-dominated diffusion problems \cite{demkowicz2013robust,ellis2016robust,chan2014robust,chan2013dpg}, and PDE-constraint optimization problems \cite{bui2014pde, causin2005discontinuous} from computational mechanics.

Designing the DPG method with optimal test space which is different from trial space, as a projection of trial space, at a continuous space, implies continuity and coercivity of the discrete scheme under adequate regularity characteristic of test and trial space on any mesh. Moreover, the automatic adaptive version of this method is guaranteed by a built-in error indicator. However, to the best of the author's knowledge, this method has not yet been used in the quantitative finance community in spite of its potential benefits, and so many more interesting characteristics of the method that one can find in the literature. 
 
In this paper, owing to the unconditional stability and solid mathematical theory of the DPG method, we proposed this method for the problem of option pricing and estimating Greeks under the Black-Scholes model. We propose both ultraweak and primal formulations of the DPG method for pricing Vanilla options, American options, Asian options, and double Barrier option and their sensitivity analysis. However, the time-space DPG is not the target here and a time-stepping strategy is used for solving the problem through time. Graph norm in which the optimal test space is established is designed and through different numerical examples, the efficiency of the proposed methods is assessed for both ultraweak and primal formulation. Besides, the DPG method for the free boundary value problem and linear complementarity (LCP) problem corresponding to the American option is provided and the early exercise boundary is obtained accordingly. Finally, the sensitivity analysis of the option price, Greeks, with respect to the underlying parameters are evaluated.  
 
 Computing the optimal test space through the test-to-trial operator introduced in the original mathematical theory of the method \cite{chan2013dpg}, \cite{roberts2013discontinuous} makes this method relatively computationally expensive. However, utilizing a broken test space overcomes this issue by localizing the evaluation of optimal test space that is conforming element-wise. Using the method with discontinuous optimal test space will allow parallelizing the assembly of the computation, and alongside local computation of test space makes the method reliable and viable. This feature of the DPG method can help to develop a high-performance implementation of the method and take advantage of highly parallel computers. Recently, some effort has been made in the form of designing a software framework to simplify the implementation of the DPG method. Camellia \cite{roberts2014camellia} is a C++ software introduced by Nathan V. Robert in the Argonne National Laboratory to allow developers to create a hp-adaptive DPG method. Astaneh et. al \cite{astaneh2018high} proposed PolyDPG in MATLAB to implement the polygonal DPG method using ultraweak formulation. However, in this paper, a prototype high-performance C++ code is developed independent of the aforementioned software for option pricing using both ultraweak and primal DPG formulation which is accessible to the public in
\href{https://github.com/DavoodDamircheli/HPCUltraweakDPGforOptionPricing}{HPC DPG for Option pricing}. 

It is worth noticing that our intention is not to compete with the previous numerical schemes used in the literature despite all the desirable aspects of the DPG method. Using this method more widely by researchers, the features of this method might be handier in more complicated and more challenging problems in quantitative finance including option pricing in higher dimensions than one-dimensions.

The outline of this paper is as follows. A very brief introduction of the Discontiouse Petrov Galerkin method with optimal test space is present in section \ref{DPGmethod}, and notation and elementary tools from the functional analysis are set in section \ref{FuncSpac}. In the sequel, we first present the DPG method for the vanilla option in section \ref{VanlOption}. We will introduce the graph norm of the DPG method for both primal and ultraweak formulation in this section. An experimental base convergent analysis is conducted for the European option pricing. Exotic option pricing including American option, Asian option, and a double barrier is numerically solved in section \ref{exotcOption}. in this part of the paper, we introduced the graph norm pertaining to every exotic option considered in this context. Standard examples in the literature are investigated with the proposed method. Finally, the DPG method is presented for sensitivity analysis of the option pricing problem in section \ref{greeks}, and the performance of the DPG method for evaluating Greeks for both exotic and vanilla options are assessed in the section. 

\newpage
%------------------------------ section I --------------------------------------
\section{The DPG Method}\label{DPGmethod}
In this section, we briefly provide a high-level introduction to the Discontinuous Petrov-Galerkin Method with Optimal Test Function. A review of the method is given for the steady-state problem, and the transient version of the method with a more concrete definition of the spaces to treat the option problem will be presented in section (\ref{VanlOption}).
Let's begin with the standard well-posed abstract variational formulation which has not necessarily symmetric functional setting, seeking $u\in U$ such that 
\begin{align}\label{stdSytem}
    b(u,v)=l(v), \quad v\in V,
\end{align}
where trial space $U$ and test space $V$ are proper Hilbert spaces. $l(\cdot)$ is a continuous linear functional, $b(\cdot, \cdot)$ is a bilinear (sesquilinear) form that satisfies the inf-sup condition as follows:
\begin{align}\label{infSupCont}
    \sup_{v\in V} \frac{|b(u,v)|}{\|v\|_V}\geq \gamma \|u\|_U, \quad \forall u\in U,
\end{align}
which guarantees the well-posedness of the variational form (\ref{stdSytem}). Therefore, discretize version of variational form (\ref{stdSytem}) with Petrov-Galerkin method is problem of finding $u_h\in U_h\subset U$ such that 
\begin{align}\label{stdPetrovGalSytem}
    b(u_h,v_h)=l(v_h), \quad v_h\in V_h.
\end{align}
Based on Babu$\check{s}$ka's theorem (\cite{babuvska1971error}) for a discretized system (\ref{stdPetrovGalSytem}) in a case where $\dim(U_h)=\dim(V_h)$, is stable or to another word the system is well-posed if the discrete inf-sup condition is satisfied as follows 
\begin{align}\label{infSupDis}
    \sup_{v_h\in V_h} \frac{|b(u_h,v_h)|}{\|v_h\|_V}\geq \gamma_h \|u_h\|_U, \quad \forall u_h\in U,
\end{align}
where the inf-sup constant $\gamma_h$ must be bounded away from zero meaning
$\gamma_h\geq\gamma >0$. Now, choosing the discrete spaces of trail and test space is of matter of importance. Indeed, trial space $U_h$ is usually picked by approximability, but trial space $V_h$ can be chosen in such a way as to dictate special properties of the numerical algorithm such as being well-posed.

The Petrov-Galerkin method with optimal test space has been designed in a way that for each discrete function $u_h$ from trial space $U_h$, it finds a corresponding optimal test function $v_h\in V$ as a supremizer of inf-sup condition, i.e optimal test function $v_h\in V$ construct such that 
\begin{align}
    \sup_{v\in V} \frac{|b(u,v)|}{\|v\|_V}= \frac{|b(u,v_h)|}{\|v_h\|_V}.
\end{align}
Given any trial space $U_h$, let's define a trial-to-test operator $T:U\longrightarrow V$. The optimal test space is defined as the image of the trail space via this operator $V_h^{\text{opt}}:=T(U_h)$, where the function from optimal test space $v^{\text{opt}}\in V_h^{\text{opt}}$ is satisfying in 
\begin{align}\label{trialToTestOpt}
    (v^{\text{opt}},v)_V=(Tu_i,v)_V=b(u_i,v), \quad \forall v\in V,
\end{align}
in which $(\cdot,\cdot)_V$ is the inner product on the test space. In fact, the equation (\ref{trialToTestOpt}) uniquely determines the optimal test space with the Riesz representation theorem with which discrete stability of the discrete form (\ref{stdPetrovGalSytem}) automatically is attained. The test function defined in (\ref{trialToTestOpt}) is designed in a way that the supremizer of the inf-sup continuous condition implies the satisfaction of the discrete inf-sup condition and as a result, it guarantees the discrete stability. Moreover, we will have
\begin{align}
    \sup_{v_h\in V_h^{\text{opt}}} \frac{|b(u_h,v_h)|}{\|v_h\|_V}\geq \frac{|b(u_h,Tu_h)|}{\|Tu_h\|_V}=\sup_{v\in V}\frac{|b(u_h,v)|}{\|v\|_V}\geq \gamma \|u_h\|_U,
\end{align}
so, we have inf-sup constant $\gamma_h\geq \gamma$.
\begin{thm}
  The trial to test operator $T:U\longrightarrow V$ is defined by:
  \begin{align}
      Tu = R_V^{-1}Bu, \quad u\in U
  \end{align}
  where $R_V:V\longrightarrow V'$ is the Riesz operator corresponding to test inner product. In particular, T is indeed linear.
\end{thm}
\begin{proof}
  see \cite{demkowicz2020oden}
\end{proof}
It can be shown (\cite{demkowicz2020oden}) that the Ideal Petrov-Galerkin method introduced above is equivalent to a mixed method as well as a minimum residual method where residual is defined in a dual norm. The ideal PG method benefits from a built-in error indicator for mesh adaptivity thanks to the corresponding mixed method where Riesz's representation of the residual in the dual test norm has been exploited. Assume $\epsilon$ is the solution of the following variational form for a given $u_h\in U_h$
\begin{align}\label{RisRpMix}
    (\epsilon,v)_V=l(v)-b(u_h,v),\quad \forall v\in V.
\end{align}
So, the Riesz representation of the residual $\epsilon$ is uniquely defined by (\ref{RisRpMix}). Then the following mixed problem can be defined 
\begin{equation}\label{mixedGal}
    \left\{\begin{split}
         &u_h\in U_h, \quad \epsilon \in V,\\
         \\
         &(\epsilon,v)_V+b(u_h,v)=l(v), \quad v\in V,\\
         \\
         &b(\delta u_h,\epsilon)= 0, \quad \delta u_h\in U_h,
    \end{split}\right.
\end{equation}
where the solution of the Ideal Petrov-Galerkin problem with optimal test space can be derived from solving the mixed Galerkin problem (\ref{mixedGal}). Thus, the method inherently has a built-in residual a-posteriori error $\epsilon$ measured in the test norm.

Nevertheless, determining the optimal test functions analytically except for some simple model problems is impossible. Therefore, to some extent approximating optimal test space in a way that the discrete inf-sup condition satisfies, is a necessity. An enriched test subspace $V_h\subset V$ is exploited as a remedy for this approximation. So, the Practical Petrov-Galerkin method with optimal test space approximated by enriched test space can be obtained as follows:
\begin{equation}
    \left\{\begin{split}
         &u^r_h\in U_h,  \\
         \\
         &b(u^r_h, T^r \delta u_h)=l(T^r \delta u_h),\quad \delta u_h\in U_h, 
    \end{split}\right.
    \end{equation}
where approximated optimal test space computes with component satisfy the standard discretization
\begin{equation}
    \left\{\begin{split}
         &T^ru\in V^r,  \\
         \\
         &( T^r u,  \delta u_h)_V = b(u, \delta v),\quad \delta v\in V_r. 
    \end{split}\right.
\end{equation}
Indeed, we increase the dimension of the discrete enriched test space in order to meet the discrete inf-sup condition for the system (\ref{stdPetrovGalSytem}). This strategy is valid due to Brezzi's theory \cite{demkowicz2020oden} that allows the dimension of discrete test space $V^r$ exceed the dimension of the trial space in spite of Babu$\check{s}$ka's theory which enforces the dimension of discrete trial and dimension of discrete test space to overlap. Analysis of stability reduction in practical Petrov-Galerkin method can be performed exploiting Fortin operators \cite{gopalakrishnan2014analysis}, \cite{nagaraj2017construction}. 

In spite of the myriad of advantages that the practical Petrov-Galerkin Method introduced so far, due to the computation of optimal test space globally through operator $T$, it is very expensive. Utilizing a broken test space overcomes the issue of localizing the evaluation of optimal test space that is conforming element-wise. Therefore, using the method with discontinuous optimal test space will parallelize the assembly of the computation alongside the local computation of test space making the method reliable and viable. Besides, this will justify the name of the Discontinuous Petrov-Galerkin method (DPG) with optimal test functions. However, breaking the test space will bring the need for introducing additional trace variables and flux variables on the mesh skeleton on the element interface. This will discuss thoroughly in section (\ref{VanlOption}) by proposing the DPG method on the Ultraweak and primal formulation of the option pricing problem. 
%--------------------------------Section2 ---------------------------------------

\section{Functional Spaces and Preliminaries}\label{FuncSpac}
we define following energy spaces to deal wit option pricing problem 
\begin{equation}\label{omegaSpace}
    \begin{split}
        &L^2(\Omega) = \{f:\Omega\to \mathbb{R}\quad|\quad\|f\|_{L_2}\leq \infty\},\\
        &H^1(\Omega) = \{f:\Omega\to \mathbb{R}\quad|\quad f\in L^2(\Omega)  \quad 
        \Delta f\in(L^2(\Omega))^d\},
    \end{split}
\end{equation}
with the $L_2$-norm defined as 
\begin{align}
    \|f\| :=(f,f)^{\frac{1}{2}}( \int_{\Omega}|f|^2dx )^{\frac{1}{2}}.
\end{align}
The domain of the problem $\Omega$ is partitioned into a set of computational domain $\Omega_h$ with open disjoint elements of $\{L\}_{L\in \Omega_h}$. This computational domain of trial space for all problems is $[-6,6]$, except for the Asian option which is $[-2,2]$. Having the finite element mesh $\Omega_h$, we can define corresponding broken energy space as
\begin{equation}\label{omegaSpace1}
    \begin{split}
        &L^2(\Omega_h) = \{f\in L^2(\Omega)\quad|\quad|f|_{L}\in L^2(L),\forall L\in \Omega_h\},\\
        &H^1(\Omega_h) = \{f\in L^2(\Omega)\quad|\quad|f|_{L}\in H^1(L),\forall L\in \Omega_h\},
    \end{split}
\end{equation}
using broken test space, we need to define the energy space for the trace variable as well. We define these spaces on the mesh skeleton $\Gamma_h$,as 
\begin{align}
    H^{\frac{1}{2}}(\Gamma_h)=\{\hat{f}\in \prod_{L\in \Omega_h}  H^{\frac{1}{2}}(\partial L)\quad | \quad \exists y\in H^{1}(\Omega)\quad \text{s.t} \quad \phi(y\big|_L) = \hat{f}           \},
\end{align}
where the operator $\phi(\cdot)$ is continuous trace operator can be defined element-wise
\begin{align}
\phi :H^1(\Omega_h)\to \prod_{L\in \Omega_h}  H^{\frac{1}{2}}(\partial L).    
\end{align}
Moreover, we need to define appropriate space for the variational inequality defined on problem of American option pricing. Thus, we define a half space $\mathcal{H}$ as following
\begin{align}
    \mathcal{H}:=\{f\in L^2(\mathbb{R}_{+})\quad | \quad f\geq f^{*}\},
\end{align}
where  $f^{*} \in L^2(\mathbb{R}_{+}) $ is the obstacle function. One can see more detail of this space in \cite{achdou2005computational}, \cite{tremolieres2011numerical}.
It is worth noticing that for option pricing in one-dimension we consider a uniform discretization of time interval $[0,T]$, truncated domain of space $[x_{\min},x_{\max}]$ as the finite element mesh $\Omega_h$.

%-----------------------------Option Pricing---------------------------------------
\section{Pricing Vanilla Options}\label{VanlOption}
In this section, we use the DPG method introduced in section \ref{DPGmethod} to numerically solve the option pricing problem. pricing vanilla option with the DPG method is presented in this section, and the exotic option is investigated in the next section.
%---------------------------------Black-Scholes---------------------------------------
\subsection{Vanilla European Options Based on Black-Scholes Model}\label{BlkSholes}
In this part, we use the DPG method for the popular Black-Scholes Model which simply provides a closed-form solution to all European-type derivatives (vanilla option). It is worth mentioning that even though assumptions of this model are not worldwide valid, there are still a large group of people on the market that will use the Black-Scholes model plus a premium \cite{higham2004introduction}. Besides, this model can be used as a test model to assess the efficiency of the method. Let's recall the Black-Scholes model briefly. This model assumes that the price of a risky asset, $S_t$, is evolving as a solution of the stochastic differential equation as follows
\begin{align}\label{geoBrownMo}
    dS_t = r S_t dt + \sigma_t  S_t dW_t,
\end{align}
in which $W_t$ is the Wiener process on a appropriate probability space $(\Omega, \mathcal{A},\mathbb{P},\mathcal{F}_t)$, $r$ is a risk free interest rate, and $\sigma_t$ is volatility of the return on the underlying security. The SDE (\ref{geoBrownMo}) is called geometric Brownian motion as well. Let's consider a European style call option on an underlying asset $S_t$, where this spot price $S_t$ satisfies in the geometric Brownian motion like (\ref{geoBrownMo}) and with the payoff of $\max \{ S_T-K,0\} = (S_T-K)_{+}$ at the date of expiration $T$ for the striking price $K$. We are interested in the fair price of this option at the current moment,$U(S_{0},0)$, if we denote the value of the option by $U(S_t,t)$ at time $t$, when the underlying price has the value $S_t$. The Black-Scholes formula express value of the option as 
\begin{align}
    U(S_t,t) = \mathbb{E}^{Q}(e^{-\int{t}^{T}r_tdt}(S_T-K)_{+}|\mathcal{F}_t),
\end{align}
It can be shown \cite{achdou2005computational,higham2004introduction} that option price of $U(S_t,t)$ satisfies in the followings deterministic partial differential equation. 
\begin{align}\label{BlS1}
    \frac{\partial U}{\partial t}+\frac{\sigma^2}{2}S_t^2\frac{\partial U^2}{\partial S^2}+r S \frac{\partial U}{\partial S}-rU(S,t) = 0,
\end{align}
with the following boundary condition 
\begin{equation}\label{BC1}
\begin{split}
&U(0,t) = 0, \quad \forall t\in [0,T],  \\
       \\
&\lim_{S_t \to \infty}U(S_t,t) = S_t-e^{-r(T-t)}, \quad \forall t\in [0,T].  
\end{split}
\end{equation}
It is well-known \cite{shreve2004stochastic,higham2004introduction,achdou2005computational}, having the upper tail of the standard normal distribution
\begin{align}
    N(x) = \frac{1}{\sqrt{2 \pi}}\int_{-\infty}^{x}e^{-\frac{z^2}{2}}dz,
\end{align}
and 
\begin{align}
    d1 &= \frac{\log(S_0/K)+(r+\frac{\sigma^2}{2}T)}{\sigma \sqrt{T}},
    &
    d2 &= d1-\sigma \sqrt{T},
\end{align}
the solution of equation (\ref{BlS1}) for a European call option can be expressed as
\begin{align}\label{BlSanalytic}
    U(S_t,t) = S_tN(d_1)-Ke^{-r(T-t)}N(d2).
\end{align}
The closed-form analytical solution (\ref{BlSanalytic}) for the European call option is used as a benchmark to study the accuracy and efficiency of the DPG method. Switching log-prices $x = \log(\frac{S_t}{S_0})$, and changing variable $\tau = T-t$, the partial differential equation (\ref{BlS1}), and the boundary conditions (\ref{BlS1}) can transfer to the following initial value constant coefficient partial differential equations solution of equation (\ref{BlS1}) for a European call option can be expressed as 
\begin{equation}\label{BlSconst}
    \left\{
    \begin{split}
        &\frac{\partial U}{\partial \tau}-\frac{\sigma^2}{2}\frac{\partial U^2}{\partial x^2}-(r +\frac{\sigma^2}{2})\frac{\partial U}{\partial x}+rU(x,t)= 0,  \\
         \\
    &U(x,0) = \max(e^x-K,0), \\
         \\
       &U(0,\tau) = 0,\\
         \end{split}
    \right.
\end{equation}
Noting that equation (\ref{BlSconst}) can be used for pricing of derivatives whose payoff depends on the price of the underlying asset at the maturity date, and more complicated options whose price are path-dependent such as American options and Asian options will use different approaches that we present them in the coming sections. 
We use finite-difference $\theta$-method to discretize the time derivative of the problem (\ref{BlSconst}) with the following form 
\begin{align}\label{TimeDisc}
    \frac{u^{n+1}-u^{n}}{\Delta \tau}-(\theta \mathcal{L}_{BS}u^{n+1}+(1-\theta)\mathcal{L}_{BS}u^{n}) = 0,
\end{align}
for $n=0,1,2,N_{\tau}-1$, with the time step $\Delta  \tau = T/N_{\tau}$, and implicitness factor $\theta\in [0,1]$.  
Besides, operator $\mathcal{L}_{BS}$ is defined as follows
\begin{align*}\label{BSoperator}
    \mathcal{L}_{BS}u = -\frac{\sigma^2}{2}\frac{\partial^2u(x,t)}{\partial x^2}-(r +\frac{\sigma^2}{2})\frac{\partial u(x,t)}{\partial x}+r u(x,t),
\end{align*}
So, different values for $\theta$ can lead us to different well-known time-stepping schemes such Backward Euler method ($\theta = 1.0$), Crank-Nicolson method ($\theta = 0.5$), and forward Euler method ($\theta = 0.0$). The Numerical efficiency of the finite difference method is well-known in the literature \cite{bulirsch2002introduction}. We proceed with introducing the DPG methodology for spatial discretization. Varieties of the variational formulation can be developed for the semi-discrete model problem (\ref{TimeDisc}) with different properties. In this investigation, we concentrate on two formulations including the classical (primal) formulation and the ultraweak formulation.
%-------------------- primal for vanilla options -------------------------------
\subsection{Primal formulation for Vanilla options}
In this subsection we propose the standard classical varational formulation for DPG method that is called the DPG primal formulation. Testing semi-discrete problem (\ref{TimeDisc}) with a proper test function $v$, integrating over the domain and using Green identity, we will have
\begin{equation}
\begin{split}
    &(u^{n+1},v)-(u^{n},v)\\
    &\qquad \Delta \tau \theta\Big[-(\frac{\sigma^2}{2}\frac{\partial}{\partial x} u^{n+1},\frac{\partial}{\partial x}v)_{\Omega_h} +((r +\frac{\sigma^2}{2})\frac{\partial}{\partial x} u^{n+1},v)_{\Omega_h}-
    (r u^{n+1},v)_{\Omega_h}+
    \langle \frac{\partial}{\partial x} u^{n+1},v\rangle_{\partial \Omega_h}\Big]\\
    &\qquad+\Delta \tau (1-\theta)\Big[-(\frac{\sigma^2}{2}\frac{\partial}{\partial x} u^{n},v)_{\Omega_h} +((r +\frac{\sigma^2}{2})\frac{\partial}{\partial x} u^{n},v)_{\Omega_h}-
    (r u^{n},v)_{\Omega_h}+\langle \frac{\partial}{\partial x} u^{n+1},, v\rangle_{\partial \Omega_h}\Big]=0,\\
\end{split}
\end{equation}
where $(\cdot,\cdot)$ are standard inner product in the Hilbert space $L_2$  and $\langle\cdot,\cdot\rangle$ is the duality pair in the $L^2(\Gamma)$. Trial space is tested with a broader discontinuous (broken) space in the DPG methodology, so as a result we don't assume that test functions disappear on the Dirichlet boundary conditions. However, the term $\frac{\partial u}{\partial x}$ will be recognized as the flux variable $\hat{q}_n$ which is a new unknown on the mesh skeleton. Thus, Defining a new group variable ${\bf u}=(u, \hat{q})\in H_1(\Omega)\times H^{-1/2}(\partial \Omega)$, the broken primal formulation for Black-Scholes (\ref{BlSconst}) reads
\begin{equation}\label{primalform}
    \left\{ 
    \begin{array}{l}
         b_{\text{primal}}({\bf u},v)=l(v),  \\
         \\
          {\bf u}(e^x,0) = \max(x-K,0),\\
          \\
          {\bf u}(0,\tau) = 0.
    \end{array}
    \right.
\end{equation}
where bilinear form $b_{\text{primal}}(\cdot,\cdot)$ and  linear operator $l(\cdot)$ are defining as follows
\begin{equation}\label{primalBS}
%\left\{
\begin{split}
&b_{\text{primal}}({\bf u},v) = (u^{n+1},v)+\Delta \tau \theta\Big[(-\frac{\sigma^2}{2}\frac{\partial}{\partial x} u^{n+1},\frac{\partial}{\partial x} v)_{\Omega_h} +((r +\frac{\sigma^2}{2})\frac{\partial}{\partial x} u^{n+1},v)_{\Omega_h}\\
&\qquad-(u^{n+1},v)_{\Omega_h}+ \langle \hat{q}^{n+1}, v\rangle_{\partial \Omega_h}\Big], \qquad \quad  n=1,\cdots, N_t,\\
\\
&l(v) = (u^{n},v)+\Delta \tau (1-\theta)\Big[(\frac{\sigma^2}{2}\frac{\partial}{\partial x} u^{n},\frac{\partial}{\partial x}v)_{\Omega_h} +((r +\frac{\sigma^2}{2})\frac{\partial}{\partial x} u^{n},v)_{\Omega_h}\\
&\qquad-(u^{n},v)_{\Omega_h}-\langle \hat{q}^{n}, v\rangle_{\partial \Omega_h}\Big],\qquad \quad  n=1,\cdots, N_t,\\
\end{split}
%\right.
\end{equation}
and boundary conditions ${\bf u}^{0} = \max(e^x-K,0), \forall x\in \Omega_h$, and ${\bf u}^i(x=0) = 0, \forall i=1,\cdots,N_t$.
Let's recall that here element-wise operations are denoted by subscribing $h$. 
Having the new flux unknown on the mesh skeleton in the primal formulation (\ref{primalBS}) is the price that we pay to use a larger test space (enriched test space).   
%--------------------------------Ultraweak formulation ------------------------
\subsection{Ultraweak Formulation for Vanilla Options}
In this section, we will derive the ultraweak DPG formulation for the pricing problem. The first step for finding ultraweak formulation is to transform the Black-Scholes problem (\ref{BlSconst}) into a first-order system of differential equation by defining a new variable $\vartheta(x,t) =\frac{\partial U}{\partial x}(x,t), \quad \forall (x,t)\in \Omega \times [0,T] $ as following 
\begin{equation}\label{firstOrderBS}
    \left\{
    \begin{split}
        &\frac{\partial U}{\partial \tau}-\frac{\sigma^2}{2}\frac{\partial  {\bf \vartheta}}{\partial x}-(r +\frac{\sigma^2}{2})\vartheta+rU(x,t)= 0,  \\
         \\
         &{\bf \vartheta}-\frac{\partial U}{\partial x}=0,\\
         \\
          &  U(x,0) = \max(e^x-K,0), \\
         \\
       &U(0,\tau) = 0.\\
         \end{split}
    \right.
\end{equation}
By defining a new group variable ${\bf u} = (u,\vartheta)$, testing the equation (\ref{firstOrderBS}) with the test variables ${\bf v} = (v,\omega)$, and integrating and using Green's identity, we will have 
\begin{equation}\label{ultraVanila}
\begin{array}{l}
    (u^{n+1},v)+(u^{n},v)+\\
    \Delta \tau \theta
    \bigg[({\bf \vartheta}^{n+1},\frac{\sigma^2}{2}  \frac{\partial}{\partial x} v)_{\Omega_h} +({\bf \vartheta}^{n+1},(r +\frac{\sigma^2}{2}) v)_{\Omega_h}-(u^{n+1},v)_{\Omega_h}+
    -(u^{n+1},\frac{\partial }{\partial x}\omega)
    -({\bf \vartheta}^{n+1},\omega)+\\
    \\
    \hspace{2cm }\langle \frac{\partial}{\partial x} u^{n+1}, v\rangle_{\partial \Omega_h}
    \langle \frac{\partial}{\partial x} {\bf \vartheta }^{n+1}, v\rangle_{\partial \Omega_h}\bigg]
   +\Delta \tau \theta
    \bigg[({\bf \vartheta}^{n},\frac{\sigma^2}{2}  \frac{\partial}{\partial x} v)_{\Omega_h} +({\bf \vartheta}^{n},(r +\frac{\sigma^2}{2}) v)_{\Omega_h}-\\
    \\
    \hspace{2cm} (u^{n},v)_{\Omega_h}+
    -(u^{n},\frac{\partial }{\partial x}\omega)
    -({\bf \vartheta}^{n},\omega)
   +\langle \frac{\partial}{\partial x} u^{n}, v\rangle_{\partial \Omega_h}
    +\langle \frac{\partial}{\partial x} {\bf \vartheta }^{n}, v\rangle_{\partial \Omega_h}\bigg] =0,
\end{array}
\end{equation}
As has been noted above we use a discontinuous test space where this space is element-wise conforming to the DPG methodology. Besides, in ultraweak formulation, there is no derivative of the trial variable in this weak formulation, and these trial variables are defined in $L_2(\Omega)$, therefore, the boundary values of the field variables are meaningless on the skeleton $\Gamma$. Thus, we introduce two trace variables $\hat{u}_{n+1}\in H^{1/2}(\Omega)$, and $\hat{\vartheta}^{n+1}\in H^{1/2}(\Omega)$ that are unknown on the skeleton. If we define the group variables ${\bf u} = (u,{\vartheta})$, ${\bf \hat{u}} = (\hat{u},\hat{{\vartheta}})$, and ${\bf v} = (v,\omega)$, the broken ultraweak formulation corresponding to the Black-Scholes model will be finding ${\bf u} = (u,\vartheta)\in L_2(\Omega)\times L_2(\Omega)$, and ${\bf \hat{u}} = (\hat{u},\hat{{\vartheta}})\in H^{1/2}(\Omega)\times H^{1/2}(\Omega) $ such that
\begin{equation}\label{ultrweakForm}
    \left\{ 
    \begin{array}{l}
         b_{\text{ultraweak}}(({\bf u},\hat{{\bf u}}),{\bf v})=l({\bf v})  \\
         \\
          ({\bf u},\hat{{\bf u}})|_{(x,0)} = \max(e^x-K,0), \\
         \\
          ({\bf u},\hat{{\bf u}})|_{(0,\tau)} = 0,\\
    \end{array}
    \right.
\end{equation}
where
\begin{align}
\begin{split}
    &b_{\text{ultraweak}}(({\bf u},\hat{{\bf u}}),{\bf v}) = 
    b_{\text{ultraweak}}(((u,{\vartheta}), (\hat{u},\hat{{\vartheta}})), (v,\omega))
    \\
    \\
    &\hspace{3cm}=(u^{n+1},v)+\Delta \tau \theta\bigg[({\bf \vartheta}^{n+1},\frac{\sigma^2}{2}  \frac{\partial}{\partial x} v)_{\Omega_h} +({\bf \vartheta}^{n+1},(r +\frac{\sigma^2}{2}) v)_{\Omega_h}-(u^{n+1},v)_{\Omega_h}-\\
    \\
   &\qquad(u^{n+1},\frac{\partial }{\partial x}\omega)
    -({\bf \vartheta}^{n+1},\omega)+
        \langle \hat{u}^{n+1}, v\rangle_{\partial \Omega_h}+
    \langle \hat{{\vartheta}}^{n+1}, v\rangle_{\partial \Omega_h}\bigg],\qquad n=1,\cdots, N_t\\
   \\
     &l({\bf v}) = l(v,\omega)= (u^{n},v)+\Delta \tau \theta
    \bigg[({\bf \vartheta}^{n},\frac{\sigma^2}{2}  \frac{\partial}{\partial x} v)_{\Omega_h} +({\bf \vartheta}^{n},(r +\frac{\sigma^2}{2}) v)_{\Omega_h}
       -(u^{n},v)_{\Omega_h}+\\
   \\
    &\hspace{2cm}-(u^{n},\frac{\partial }{\partial x}\omega)
    -({\bf \vartheta}^{n},\omega)+\langle \hat{u}^{n}, v\rangle_{\partial \Omega_h}
    +\langle \hat{{\vartheta}}^{n}, v\rangle_{\partial \Omega_h}\bigg], \qquad n=1,\cdots, N_t,
\end{split} 
\end{align}
with the boundary condition $({\bf u},\hat{{\bf u}})^{0} = \max(x-K,0), $, and $({\bf u},\hat{{\bf u}})^{i} = 0, n=1,\cdots, N_t$.
It is well-known fact \cite{demkowicz2010class,demkowicz2011analysis}, that the DPG method significantly depends on the choice of the test space's inner product since it determines the norm and as a result the structure of test space in which the DPG method is optimal. As an illustration, if the errors in $L_2$-norm are of interest, there is a tangible theory \cite{demkowicz2020oden} that shows that the graph norm is a suitable choice for the test space in ultraweak formulation, and the standard energy norm induced form bilinear $\| \cdot\|_E = b_{\text{primal}}(v,v)$ is the candidate the primal formulation. we employ the following test norms for the formulations proposed above. 
In this paper, we propose the following graph norm  (\ref{ultrweakForm}), and (\ref{primalform})
\begin{align}\label{testNormVanial}
    \begin{split}
         &\text{Primal}: \|v\|^2_{V} = \frac{1}{\Delta t} \|v\|^2+\frac{1}{(\Delta t)^2} \|\frac{\sigma^2}{2}\frac{\partial }{\partial x}v\|^2,  \\
         \\
         &\text{Ultraweak}: \|{\bf v}\|^2_{V} = \|(v,\omega)\|^2_{V}=\frac{1}{(\Delta t)^2} \| \frac{\sigma^2}{2}\frac{\partial }{\partial x} v- r v-\omega\|^2+ \frac{1}{\Delta t}\| (r +\frac{\sigma^2}{2}) v-\frac{\partial }{\partial x}\omega\|^2,  \\ 
    \end{split}
\end{align}
Having the graph norm and energy norm defined in (\ref{testNormVanial}), and the inner product of the corresponding test space as a direct result of it, we are ready to discretize the weak forms and construct the DPG system. 
In the classical Galerkin method, the convention is to choose the same discrete space for both trial and test spaces, so a squared linear system is expected. However, in the DPG method, discrete trial $U_h\subset U$ and test space $V_h \subset V$ have different dimensions. The practical DPG method with optimal test space benefits from enriched test space, meaning that $\dim V_h\geq \dim U_h$. We assume that $\{u_j\}^N_{j=1}$ , and $\{v_j\}^M_{j=1}$ are the bases of trial and test spaces respectively where $M\geq N$. In the DPG methodology, each trial space basis function $u_i$ and corresponding optimal test function $v^{\text{opt}}_i$ satisfy in the following system
\begin{align}
    (v^{\text{opt}}_i, \delta v)_V = b(u_i, \delta v), \quad \forall \delta v\in V.
\end{align}
Now let's define $M\times M$ Gram matrix $G= (G_{ij})_{M\times M}$ as 
\begin{align*}
    G_{ij} = (v_i,v_j)_V,
\end{align*}
and $N \times M$ stiffness matrix $B = (B_{ij})_{N\times M}$
\begin{align*}
    B_{ij} = b(u_i, v_j),
\end{align*}
for primal formulation finding matrix $B$ is straightforward from the bilinear form and test norm, however, calculating this matrix for ultraweak formulation can be confusing, where $B$ has the following structure 
\begin{align}
    B = \begin{bmatrix}
B_{uv} &  B_{\vartheta v} & B_{\hat{u}v} & B_{\hat{\vartheta}v}\\
B_{u{\omega}} &  B_{\vartheta {\omega}} & B_{\hat{u}{\omega}} & B_{\hat{\vartheta}{\omega}}\\
\end{bmatrix}_{N\times M},
\end{align}
and $l$ the mass matrix $l(v) = (f,r)$. We use high-order Lagrange basis of different orders to expand the trial space with order $P$, and enriched test space with order $p+\Delta p$ for $\Delta p = 2$. Thus 
The global assembly will have the following form
\begin{align}
    B^{\text{n-op}}{\bf u}_h = B^{T}G^{-1} B {\bf u}_h= B^{T} G^{-1}l = l^{\text{n-op}},
\end{align}
where discrete operators $B^{\text{n-op}}$, and $l^{\text{n-op}}$ are near-optimal mass and stiffness matrix for the DPG formula. It is worth noting that thanks to the broken structure of the test space, evaluating optimal test functions in the Gram matrix and its inversion are localized and therefore the global assembly can be paralleled, which makes the DPG method a practical method to solve the option pricing problem.
%-----------------------------numerical vanilla-----------------------------------
%-----------------------------------------------------------------------%
%------------------ price surface primal and ultraweak------------------
%-----------------------------------------------------------------------%
 \begin{figure}[!htb]
\includegraphics[width=\linewidth]{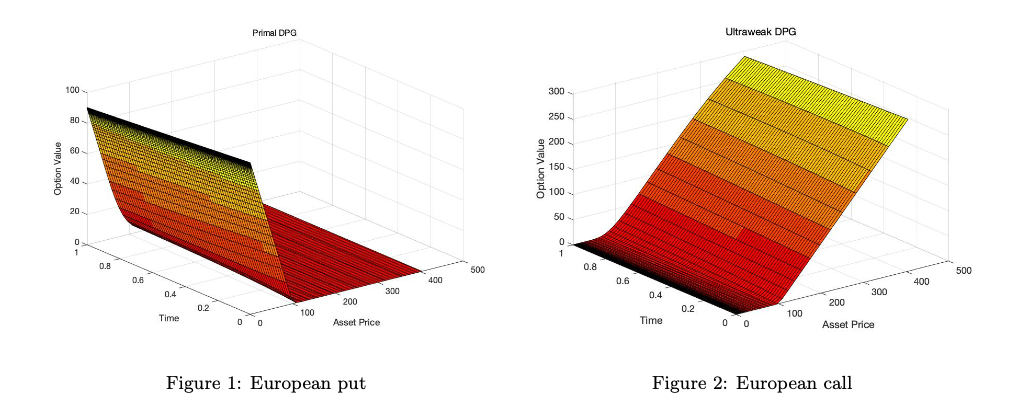}
\caption{The  surface for of two European options using DPG method with $\sigma = 0.4$, $r = 0.1$, and $k=100$. }
\label{UltraPrimaPriceSurfaceVanailla}
\end{figure}
%------------------ END price surface primal and ultraweak--------------
\subsection{Numerical Results}
In this section, we provide numerical experiments to showcase the efficiency and accuracy of the DPG method in pricing vanilla options using both the primal and ultraweak DPG methods. For this experiment, risk-free rate $r$ is set to be $0.05$, time to maturity $T$ is one year, and the strike price $K$ is $100$. The computational domain is $[-6,6]$, and a variety of values for the market volatility $\sigma$ is considered in this part. 

Through this paper, we report the relative errors of $L_2$-error, $L_{\infty}$-error of the solution obtained by the proposed numerical scheme. The binomial method implemented in \cite{higham2002nine} is utilized as a benchmark and analytical solution to compare with the approximated solution obtained with the DPG method. The relative errors are defined as follows
%-----------------------------------------------------------------------%
\begin{align}\label{errorDef}
\|E\|^2_{L_2} &= \| \frac{u-\Tilde{u}}{u}\|^2_{L_2},&\|E\|_{\infty} &= \| \frac{u-\Tilde{u}}{u}\|_{L_{\infty}},    
   \end{align}
where $\Tilde{u}$ represents the estimated value attained from the numerical method. Fig. \ref{UltraPrimaPriceSurfaceVanailla} depicts the surface of a call option with volatility $\sigma=0.4$ for both primal and ultraweak DPG formulation. 
%%%%%%%%%%%%%%%%%%%%%%%%%%%%%%%%%%%%%%%%%%%%%%%%%%%%%%%%%%%%%%%%%%%%%%%%
%-----------------order Space and Time Primal DPG European options------
%%%%%%%%%%%%%%%%%%%%%%%%%%%%%%%%%%%%%%%%%%%%%%%%%%%%%%%%%%%%%%%%%%%%%%%%
\begin{figure}[ht]
\includegraphics[width=\linewidth]{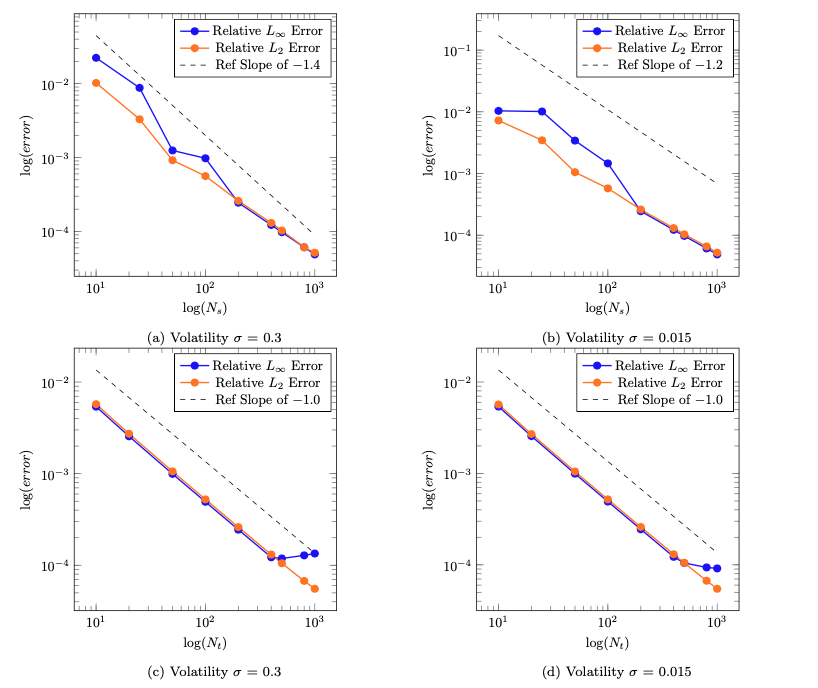}
\caption{Accuracy properties of primal DPG for European put options $r=0.05$,  $K=100$, and different volatility}
\label{orderSpacePrimalDPGvanialPut}
\end{figure}
%%%%%%%%%%%%%%%%%%%%%%%%%%%%%%%%%%%%%%%%%%%%%%%%%%%%%%%%%%%%%%%%%%%%%%%% 
%--------------------- order Space Ultraweak DPG -----------------------
%%%%%%%%%%%%%%%%%%%%%%%%%%%%%%%%%%%%%%%%%%%%%%%%%%%%%%%%%%%%%%%%%%%%%%%%
In this part of the experiment, we study the asymptotic convergence of relative errors of the numerical method for uniform mesh refinement both in time and steps. 
It is worth mentioning that error is small in general, and the relative error is of order of $10^{-6}$. 

In this regard, Fig. (2a), and (2b) displays the space order of convergence of the primal DPG method for volatilises of $\sigma = 0.3$ and $\sigma = 0.015$ pricing a European put option. It is evident that the convergence rate of primal DPG scheme is super linear in space. 

The same investigation for ultraweak DPG scheme Fig. (2a), and (2b) shows that although the convergence rate in space is super-linear the errors in this scheme decay moderately gently. We observe that for the space order in both ultraweak and primal schemes initially we see some inconsistency in the linear decreasing of the error but once a number of elements approach a certain point, we witness the expected linear convergence $\mathcal{O}(h)$, which can cause this overall super-linear convergence rate.
%-----------------------------------------------------------------------
%--------------------- order Space and Time Ultraweak DPG --------------
%-----------------------------------------------------------------------
\begin{figure}[ht]
\includegraphics[width=\linewidth]{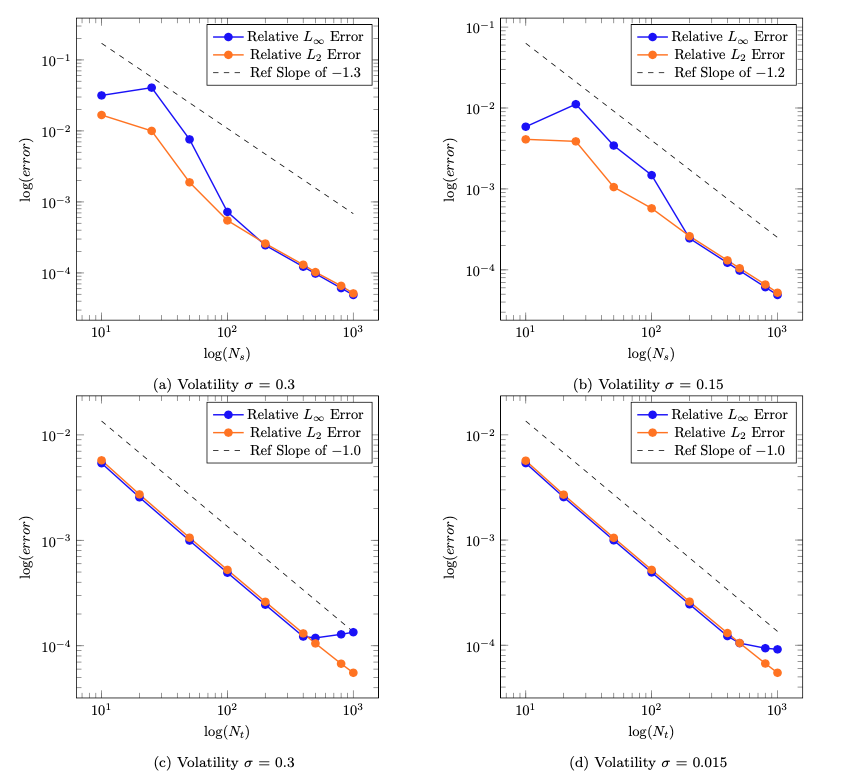}
\caption{Accuracy properties of ultraweak DPG for European put options r=0.05,  K=100, and different volatility}
\label{orderSpaceUltraweakDPGvanialPut}
\end{figure}
%------------------------End Figure---------------------------------------------
However, Fig. (3a), and (3b), and Fig. (3c), and Fig. (3d) depicts this observation more precisely when for the same scenario the rate of convergence for the Primal DPG and Ultraweak DPG method is linear in time due to the fact that the $h=0.01$ is fixed for this experiment.
%-----------------------------------------------------------------------%
%-----------------------Exotic options---------------------------------
%-----------------------------------------------------------------------%
\section{Exotic Options}\label{exotcOption}
Financial institutions issue other forms of options that are not vanilla calls or put introduced in section \ref{VanlOption}. This modern financial instrument is traded between companies and banks to cope with a variety of demands \cite{zhu2004derivative}. So, exotic options are traded in the over-the-counter (OTC) market to satisfy special needs. Being a complicated financial instrument is the common property of exotic options where the value of of these options depends on the whole or part of the path of the underlying security. Thus, exotic options are path-dependent options. In this section, we proposed the DPG method for the numerical solution of the important examples of path-dependent exotic options including American options, Asian options, Barrier options, and look-back options.
%-------------------------------American Options------------------------
\subsection{American options}\label{AmericanOptionsec}
In this section, we briefly review American option pricing under the simple model of Black-Scholes. Contrary to the European option, the holder of this contract has the right to exercise the option at any time before maturity.  
It is well known that this slight difference brings the analysis of American options much more complicated. Indeed, this right turn problem of valuing the American option into a stochastic optimization problem. The price of an American option under the risk-neutral pricing principle can be obtained as 
\begin{align}
    U(x,t) = \sup_{t\leq \tau \leq T}\mathbb{E}[e^{-\int_{t}^{\tau}r(s)ds}h(x)|\mathcal{F}_t]
\end{align}
where $h(x)$ is the option payoff, and $\tau$ is a stopping time. Stopping time is the time that owner of the option exercises the contract, besides, the stopping time is a concept in the stochastic analysis as well. It is worth noting that due to the complexity of the American option problem, this problem does not have a closed-form solution. One way of formulating American options thanks to the no-arbitrage principle is the free boundary value problem. Indeed, the free boundary happens when the option is deep in-the-money, and finding this boundary alongside pricing the American option brings extra difficulties to the problem. Here we briefly recall the different forms of American options and the corresponding DPG formulation for the formulations, for more detail one can see \cite{seydel2006tools}.

Considering the log-prices $x=\log(\frac{S_t}{S_0})$, changing tenor $T-t$ to $\tau$, the free boundary formulation of the American put option yields:
\begin{equation}\label{AmerFree}
    \left\{
    \begin{split}
       &\frac{\partial U}{\partial \tau}(x,\tau)-\frac{\sigma^2}{2}\frac{\partial U^2}{\partial x^2}(x,\tau)-(r +\frac{\sigma^2}{2})\frac{\partial U}{\partial x}(x,\tau)+rU(x,\tau) = 0,\qquad \forall x>S_f,\\
       \\
         &U(x,\tau) = K-e^x,\qquad \forall x\leq S_f,\\
         \\
         &U(x,0) = (K-e^x)^{+},\\
         \\
         &\lim_{x\to \infty}U(x,\tau) = 0,\\
         \\
         &\lim_{x\to S_f}U(x,\tau) = K-e^{S_f},\\
         \\
         &\lim_{x\to S_f}\frac{\partial U(x,\tau)}{\partial x} = -1,\\
    \end{split}
        \right.
\end{equation}
in which, $S_f$ is the free boundary of the American option pricing. It is evident that solving the problem of American option in a free boundary framework needs evaluating the free boundary along the finding the solution. Whereas, There is another approach to deriving the American option pricing problem called a linear complementarity problem (LCP). The advantage of this approach is that free boundary is not present in the formulation anymore. However, solving the LCP problem has its own complexity, and techniques  \cite{murty1988linear} . The complementarity problem of the American option can be written as
\begin{equation}\label{LCP0}
\left\{\begin{split}
&\left(\frac{\partial U}{\partial \tau}(x,\tau)-\frac{\sigma^2}{2}\frac{\partial U^2}{\partial x^2}(x,\tau)-(r +\frac{\sigma^2}{2})\frac{\partial U}{\partial x}(x,\tau)+rU(x,\tau)
\right)(U(x,\tau)-h(x))= 0,\\
\\
&    \frac{\partial U}{\partial \tau}(x,\tau)-\frac{\sigma^2}{2}\frac{\partial U^2}{\partial x^2}(x,\tau)-(r +\frac{\sigma^2}{2})\frac{\partial U}{\partial x}(x,\tau)+rU(x,\tau)
\geq 0,\\
  \\
 &   U(x,\tau)-h(x)\geq 0,\\
    \\
  &  U(x,0) = (K-e^x)^{+}.
 \end{split}
\right.
\end{equation}
The main approach here is to utilize the DPG formulation for the governing equations of the equ. \ref{LCP0}, and \ref{AmerFree}  and then consider the free boundary condition for them. The using DPG method for a (LCP) is examined before in \cite{fuhrer2018dpg} for using DPG formulation for the Signorini type problem as a contact problem. However, Thomas Fuhrer et al. in \cite{fuhrer2018dpg} proposed the ultraweak formulation of the corresponding problem, here we derive both ultraweak and primal formulation of the DPG method for the problem of American option pricing as a special case of obstacle problem. 

Now, for the DPG formulation in LCP framework, we multiply the second inequality condition in the equ. \ref{LCP0} with the smooth no-negative test functions $v\in V$ where test space is a broken convex cone and  following the same process of defining trail and flux variable presented in the  section \ref{VanlOption}, and  after some integration by part we obtain 
\begin{align}\label{VI1}
    \begin{split}
         \frac{d}{d\tau}({\bf u}, {\bf v})+ b^{\tau}({\bf u},{\bf v})\geq 0,  \\
    \end{split}
\end{align}
where bilinear form for primal formulation defies as 
\begin{align}\label{primalVI}
%\left\{
\begin{split}
    b^{\tau}_{\text{primal}}({\bf u},v) =  (-\frac{\sigma^2}{2}\frac{\partial u}{\partial x} ,\frac{\partial v}{\partial x} )_{\Omega_{+}} +\bigl((r +\frac{\sigma^2}{2})\frac{\partial u}{\partial x} ,v \bigr)_{\Omega_{+}}
    -(u,v)_{\Omega_{+}}+ \langle \hat{q}, v\rangle_{\partial \Omega_{+}},\\
 \end{split}
%\right.
\end{align}
where $\Omega_{+}$ shows the non-negative part of the domain, with a set of trial and flux variables ${\bf u}=(u,\hat{q})\in H^1(\Omega)\times H^{\frac{1}{2}}(\Omega)$ , and test variable ${\bf v} = v \in H^1(\Omega)$.
Moreover, defining trail variables ${\bf u} = (u,\vartheta)\in L_2(\Omega)\times L_2(\Omega)$, and flux variables ${\bf \hat{u}} = (\hat{u},\hat{{\vartheta}})\in  H^{1/2}(\Omega)\times  H^{1/2}(\Omega) $, one can define the bilinear form in \ref{VI1} for the ultraweak formulation as following 
\begin{align}\label{ultraWkVI}
\begin{split}
&b^{\tau}_{\text{ultraweak}}(({\bf u},\hat{{\bf u}}),{\bf v}) = 
b^{\tau}_{\text{ultraweak}}\left( \bigl((u,{\vartheta}), (\hat{u},\hat{{\vartheta}})\bigl), (v,\omega)\right),
\\
\\
&\hspace{3cm}=({\bf \vartheta},\frac{\sigma^2}{2}  \frac{\partial v}{\partial x} )_{\Omega_{+}} +({\bf \vartheta},(r +\frac{\sigma^2}{2}) v)_{\Omega_{+}}-(u,v)_{\Omega_{+}}-
(u,\frac{\partial \omega}{\partial x})_{\Omega_{+}}
\\
\\
&\hspace{3cm}-({\bf \vartheta},\omega)_{\Omega_{+}}+
\langle \hat{u}, v\rangle_{\partial \Omega_{+}}+
\langle \hat{{\vartheta}}, v \rangle_{\partial \Omega_{+}}.
\end{split}
\end{align}
It is well-known that the two variational inequality proposed with the bilinear forms \ref{primalVI}, and \ref{ultraWkVI} are the parabolic variational inequalities of the first kind that admit a unique solution \cite{kinderlehrer2000introduction}. Having well-posed variational inequality of (\ref{VI1}), we can approximate the problem in a finite-dimensional space. Thus, similar to estimating the price of vanilla options, we consider the time partition $0\leq \cdots \leq T$ of the time interval $[0, T]$, and discrete trial space $U_h\subset U$, and enriched test space $V_h\subset V$ (dim $V_h\geq$ dim $U_h$) and the corresponding basis spanned $\{u_j\}^N_{j=1}$, and $\{v_j\}^M_{j=1}$ for the aforementioned spaces. We use the backward finite difference Euler method to approximate the time derivative, and as a result, the discrete DPG for variational inequalities arising from the American option pricing problem yields
\begin{align}\label{discreteVI}
    \begin{split}
        (u^{n+1}-u^n, {\bf v})+ \Delta \tau b^{\tau}_n(u^n,{\bf v})\geq 0, \qquad \forall {\bf v} \in V_h. \\
    \end{split}
\end{align}
However, writing the $\theta$-method for the second term in left hand side of the discrete variational inequality (\ref{discreteVI}) will be performed very similarly to what is proposed for vanilla options. Let $B$ and $G$ be the stiffness and Gram matrices defined by
\begin{align}
    B_{ij}&=b_n^{\tau}(u_i,{\bf v}_j), & G_{ij}&=({\bf v}_i,{\bf v}_j)_v ,  &l_i &= (u_i,v),
\end{align}
where $(\cdot,\cdot)_v$ inner product of test space obtained from the energy norm for primal DPG and graph norm for ultraweak form introduced in (\ref{testNormVanial}). So, the discrete variational inequality (\ref{discreteVI}) is equivalent to 
\begin{equation}\label{discreteVI1}
    \left\{\begin{split}
        &B^TG^{-1}l(u^{n+1}-u^n)+ \Delta \tau B^TG^{-1}B u^{n}\geq 0, \\
        \\
        &u^{n}\geq h(x),\\
        \\
        &(u^{n}-h(x))\bigl(B^TG^{-1}l(u^{n+1}-u^n)+ \Delta \tau B^TG^{-1}B u^{n}\bigr)=0,
        \\
    \end{split}\right.
\end{equation}
for $n=1,\cdots,N_{\tau}$. Setting near the optimal discrete operators of $B^{\text{n-op}}=B^T G^{-1}B$,$l^{\text{n-op}}=B^T G^{-1}l$ discrete LCP (\ref{discreteVI1}) will attain the following form
\begin{equation}\label{discreteVI2}
    \left\{\begin{split}
    &l^{\text{n-op}}(u^{n+1}-u^n)+ \Delta \tau B^{\text{n-op}} u^{n}\geq 0,& \\
        \\
        &u^{n}-h(x)\geq 0 ,&\forall n=1,\cdots,N_{\tau}\\
        \\
        &(u^{n}-h(x))\big (l^{\text{n-op}}(u^{n+1}-u^n)+ \Delta \tau B^{\text{n-op}} u^{n}\big)=0.&
        \\
    \end{split}\right.
\end{equation}
There are different approaches to solve the discrete variational inequality (\ref{discreteVI2}) including fix-point approach, penalization method, iterative method to just name few \cite{damircheli2019solution}. 
To close the section we will present the DPG formulation for solving the free boundary value problem \ref{AmerFree}. Similar to the procedure for governing equation of vanilla options, one can test the governing equation \ref{AmerFree} with the appropriate test functions, and define the following system
\begin{equation}\label{DPGfree}
    \begin{split}
         \frac{d}{d\tau}({\bf u}, {\bf v})+ b^{\tau}({\bf u},{\bf v})=0,\qquad \forall x>S_f,   \\
         \end{split}
\end{equation}
Where the \textit{bilinear form} in the equation \ref{DPGfree} has the form of \ref{primalVI} for the primal formulation and \ref{ultraWkVI} for the ultraweak formulation. Like our approach so far, we use the Backward Euler method for time derivative and trial and test space defined for the LCP form to find the following discreet system of equation
\begin{align}\label{discreteVI?}
    \begin{split}
        (u^{n+1}-u^n, {\bf v})+ \Delta \tau b^{\tau}_n(u^n,{\bf v})= 0, \quad \forall x_h>S_f,\qquad \forall {\bf v} \in V_h. \\
    \end{split}
\end{align}
Having enough fine time discretization in the above form, using the information with one time step lag can attain a good approximation of the solution of the American option. In another word, one need to notice that the final price of the American option will find from the following implicitly boundary condition 
\begin{equation}\label{freeBound}
    u^n = \left\{\begin{array}{ll}
        \max\{h(x),u^{n-1}\}, &\forall x\in \Omega^{\mathrm{o}},  \\
         \\
         h(x), &  x=\inf{ \partial\Omega},\\
         \\
         0, & x=\sup{\partial\Omega}.
          \end{array}\right.
\end{equation}
in which $h(x)$ is the payoff of American option, Boundary conditions presented in \ref{freeBound} are necessary boundary conditions of the Problem of Valuing American option pricing. 
%--------------------------------Numerical American --------------------
%------------------LINK FOR MAGNIFING THE PLOT IN LATEX---------
%https://pgfplots.net/spy-plot/
%%%%%%%%%%%%%%%%%%%%%%%%%%%%%%%%%%%%%%%%%%%%%%%%%%%%%%%
%-------------------price in maturity -----------------
%%%%%%%%%%%%%%%%%%%%%%%%%%%%%%%%%%%%%%%%%%%%%%%%%%%%%%%
\begin{figure}[!ht]
\includegraphics[width=\linewidth]{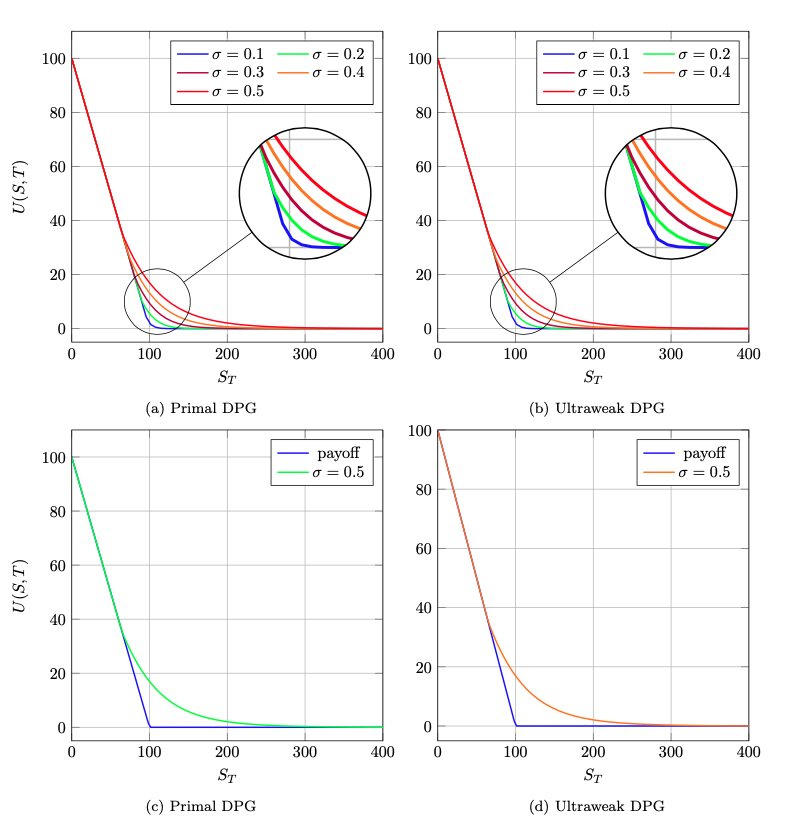}
\caption{Value of American option for r=0.05, K=100, and different volatilises}
\label{DPGdifferentVolatilityAmeicanOption}
\end{figure}
%---------------------------End of the plots---------------------------------------
\subsection{Numerical Experiments}
In this set of numerical experiments, we study the problem of valuing the American option with the ultraweak and primal DPG method. we intend to verify that DPG is a reliable and efficient method for solving this free boundary value problem. Fig. (4a), and Fig. (4b) illustrate the price of an American put option for a fixed interest rate $r= 0.05$, maturity $K= 100$, and different volatilises. It is a well-known fact that the price of an American option is greater than a European option due to the right of the owner of the American option for exercising the financial contract anytime before maturity, this can vividly be seen in Fig. (4c), and Fig. 4d for the payoff and value of an American option. Thus, the proposed methods can mimic this behavior accurately for different volatility of the market for both primal and ultraweak formulations.
%-----------------------table American Option------------------------------------
\begin{table}[htbp]
{\footnotesize
  \caption{{Value of American Option r = 0.05, $\sigma$ = 0.15, k=100 }\label{TableAmericanerror}}
\begin{center}
  \begin{tabular}{cccccccc} \hline
   $\Delta{\tau}$ & $h$  & \multicolumn{2}{c}{value} && \multicolumn{2}{c}{$\|E\|_{\infty}$}&\\ \hline
              &        &  Primal      & Ultraweak        && Primal       & Ultraweak    &\\
                              \cline{3-4}                           \cline{6-7}  
       0.01 & 0.46   &   14.           & 15.           &&     0.0159930    & 0.00963025    &\\
       0.01 & 0.23   &   14.           & 15.           &&     0.00379765  & 0.00050253    &\\
       0.01 & 0.11   &   14.           & 15.           &&     0.00074574  & 0.00133183    &\\
       0.01 & 0.05   &   14.           & 15.           &&     0.00034449  & 0.00027304    &\\
       0.01 & 0.03   &   14.           & 15.           &&     5.83E-05    & 6.21E-05    &\\
       0.01 & 0.02   &   14.           & 15.           &&     4.12E-05    & 4.48E-05    &\\
       0.01 & 0.01   &   14.           & 15.           &&     1.77E-05    & 1.75E-05    &\\ \hline
  \end{tabular}
\end{center}
}
\end{table}
%--------------------End table American Option------------------------------------
%-----------------------------end of plot---------------------------------------------
Error analysis of the American option conducted with the relative $L_2$, and $L_{\infty}$ error of the solution very similar to the definitions (\ref{errorDef}). Besides, the bench mark for the exact solution is opt the value of binomial method introduced and implemented in \cite{higham2002nine}. Table \ref{TableAmericanerror} is prepared to show the error of the DPG numerical scheme for both primal and Ultraweak formulation. In this study, the time step is fixed $\Delta \tau = 0.01$, and we use a finer mesh in spatial dimension on each step. One can see that the trend of error is descending as $h$ decreases and we get more accurate value of the American options. 
%%%%%%%%%%%%%%%%%%%%%%%%%%%%%%%%%%%%%%%%%%%%%%%%%%%%%%%%%%%%%%%%%%%%%%%% 
%------------- Order Space primal DPG for American ---------------------
%%%%%%%%%%%%%%%%%%%%%%%%%%%%%%%%%%%%%%%%%%%%%%%%%%%%%%%%%%%%%%%%%%%%%%%%
\begin{figure}[h]
\includegraphics[width=\linewidth]{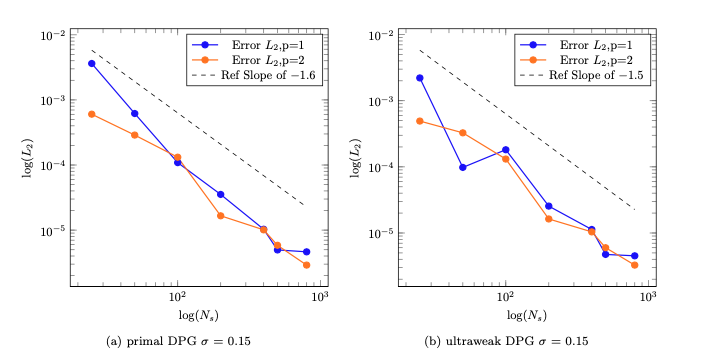}
\caption{Accuracy properties of ultraweak and primal DPG for American put options in the spatial dimension  with the parameters $r=0.05$,  $K=100$, and $\sigma = 0.15$}
\end{figure}
Although the magnitude of error is important, the order by which error is lessened is of a great importance in our error analysis. In this investigation we used the high order DPG method as well to study the effect of the order of interpolation on the valuing of the American option pricing. Let's commence with the spatial order of convergency. Fig.(5a), and Fig.(5b) illustrates the order of convergence of both primal and  ultraweak formulation for valuing American option for the fixed interest rate $r=0.05$, exercise prices of $K=100$, and the market volatility of $\sigma = 0.15$ in space order for first order and second order DPG. The experiment shows that asymptotic convergence of $L_2$ error is super linear, but it doesn't reach the $o(h^2)$ for the second order DPG scheme. One possible explanation of the diminishing the order could be an adverse impact of free boundary in the pricing problem. However, the error is relatively small, and table (\ref{TableAmericanerror}) reinforce this trend as well for relative sup-error for both primal and ultraweak formulation, where ultraweak formulation has a tiny better performance in majority of cases. 
\begin{figure}[htbp]
\includegraphics[width=\linewidth]{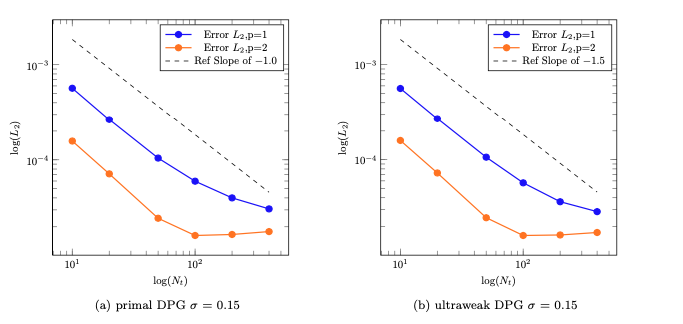}
\caption{Accuracy properties of Ultreawek and primal DPG scheme for American put options with respect to time step with parameters r=0.05,  K=100, and $\sigma = 0.15$}
%\caption{(space order )Accuracy properties of ultraweak and primal DPG for American put options r=0.05,  K=100, and \sigma = 0.15}
\label{orderTimePrimalUltraweakDPGAmericanPut}
\end{figure}
In order to study the stability and convergence in time stepping scheme, we prepared fig (\ref{orderTimePrimalUltraweakDPGAmericanPut}).A fixed mesh in space with $N_s = 64$ elements is used and decrease the time step $\Delta \tau$ and record the $L_2$-error for first and second order DPG method. The convergence analysis shows that this both primal (Fig. 6a) and ultraweak (Fig. 6b) formulation benefit from the rate of convergence of $\mathcal{O}(\Delta \tau)$ as we expected and the backward Euler method is unconditionally stable. However, the rate of convergence for time stepping captures for initial time steps (almost $N_{\tau} = 100$), where as after this point spacial discretization dictates it's impact afterwards for both DPG forms. 
%%%%%%%%%%%%%%%%%%%%%%%%%%%%%%%%%%%%%%%%%%%%%%%%%%%%%%%
%-------------------optimal exercise boundary --------
%%%%%%%%%%%%%%%%%%%%%%%%%%%%%%%%%%%%%%%%%%%%%%%%%%%%%%%
\begin{figure}[H]
\includegraphics[width=\linewidth]{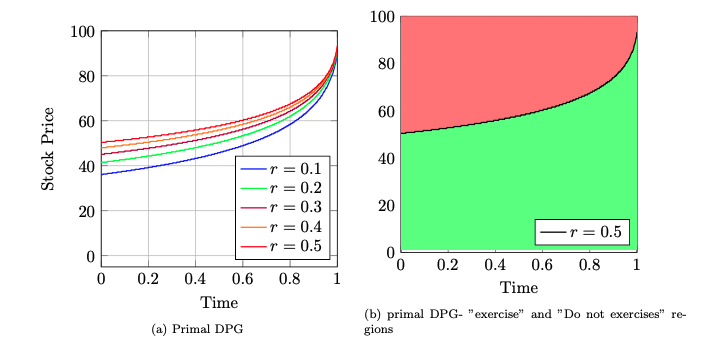}
\caption{Optimal Exercise boundary for an American put. Computed via the primal DPG method.
b) the green part is for exercise and red for "Do not exercise"}
\label{FreeBoundaryDPG}
\end{figure}
Besides accurately pricing the American-type financial derivative, finding the optimal exercise boundary for an American option is essential. The DPG method proposed in this section can find the optimal exercise boundary implicitly thanks to the projection-based method just by checking the price with the payoff at each moment or through an automatic procedure in the first active points at each time step in the primal-dual active set strategy. Fig. (7a) depicts finding this free boundary for the different interest rates of the market at each time to maturity. This optimal boundary is a powerful indicator for practitioners to choose the appropriate positions due to the hedging strategy.
Thus, the optimal exercise boundary partitions the domain of the problem into an "Exercise region" and "Do not Exercise" region (7b) where the owner of the option will exercise the option when the stock price is at the green region, and will await in the red region.

%----------------------End of Optimal exercise boundary---

%------------------------------Asian Option--------------------------------
\subsection{Asian Options}
Asian options can be classified as path-depended financial derivatives where the payoff of the option depends on the time average of the underlying security over some period of time such as the lifetime of an option \cite{shreve2004stochastic, kemna1990pricing}. This average can be taken over continuous sampling or discrete sampling and the type of average can be an arithmetic average or geometric average. The closed-form value of an Asian option is not in hand, so a numerical scheme is an essential remedy to find the value of an Asian option.

Seeking a closed-form solution such as the Laplace transform of the price for this path-dependent derivatives has been a popular approach \cite{levypricing}, \cite{vorst1992prices}, and \cite{turnbull1991quick}. However, the numerical implementation of the aforementioned methods is troublesome for low volatility cases \cite{fu1999pricing}. The Monte Carlo method can be used for the numerical solution, where it is well-known that this method is computationally expensive \cite{kemna1990pricing}, and \cite{broadie1996estimating}. Another popular approach is solving two dimensions in space PDE to find the value of an Asian option \cite{ingersoll1987theory}, \cite{vecer2001new}, and \cite{kim2014option}. Besides, Rogers and Shi \cite{rogers1995value} proposed a reduction approach where solving one-dimensional PDE obtains the value of the desired Asian option.
However, both one and two-dimensional PDEs are susceptible to oscillatory solution and can blow up through time due to existing small diffusion terms.

In this section, we propose the DPG method for pricing the option based on the Black-Scholes pricing framework. Assume the dynamic of the underlying asset satisfies in a geometries Brownian motion defined in \ref{geoBrownMo}, then the payoff of an Asian call option at maturity with the fixed-strike is following
\begin{align}\label{AsianPayoff}
    U(T) = \max\{\frac{1}{T}\int_{0}^{T}S(t)dt-K,0\}=\big(\frac{1}{T}\int_{0}^{T}S(t)dt\big)^{+}
\end{align}
based on the risk-neutral pricing theory, the price $U(t)$ of this Asian option at time $t\in[0,T]$  yields 
\begin{align}\label{expAsian}
    U(t) = \mathbb{E}[e^{-r(T-t)}U(T) | \mathcal{F}_t], \qquad \forall t\in[0,T],
\end{align}
where expectation in \ref{expAsian} is a conditional expectation with respect to the filter $\mathcal{F}_t$ of the probability space $(\Omega, P, \mathcal{F})$. Since the payoff defined in \ref{AsianPayoff} depends on the whole path of stock price $S(t)$, the price of this option is a function of $t$, $S(t)$, and the evolution of value underlying security over the path. Thus, we extend the pricing model presented in previous sections for the European and American options by defining a second process
\begin{align}
    Y(t)= \int_{0}^{t}S(v)dv,
\end{align}
where the dynamic of this new process $Y(t)$ follows a stochastic differential equation as following
\begin{align}
    dY(t)= S(t)dt.
\end{align}
Therefore, the value of the Asian option is also a function of $Y(t)$, so we denote the price of the Asian option with $U(t, S_t, Y(t))$. This function satisfies $\forall t\in[0,T]$, and $\forall (x, y)\in \mathbb{R}^{+}\times\mathbb{R}$ in the following two-dimension in space, partial differential equation(see \cite{shreve2004stochastic,kemna1990pricing} for details)
\begin{equation}
    \left\{
    \begin{split}
       &\frac{\partial U(t,S,y)}{\partial t}+\frac{\sigma^2}{2}S^2\frac{\partial U(t,S,y)^2}{\partial S^2}+r S\frac{\partial U(t,S,y)}{\partial S}
            +S\frac{\partial U(t,S,y)}{\partial y}-rU(t,S,y)=0,\\
      \\
         &U(t,0,y) = e^{-r(T-t)}(\frac{y}{T}-K)^{+},\qquad t\in[0,T],\quad y\in \mathbb{R},\\
         \\
         &U(T,S,y) = (\frac{y}{T}-K)^{+},\quad S\geq 0,\quad y\in \mathbb{R},\\
         \\
         &\lim_{y\to -\infty}U(t,S,y) = 0,\quad t\in[0,T],\quad S\geq 0.\\
           \end{split}
        \right.
\end{equation}
Now, let's define a new state variable 
\begin{align}
    x = \frac{1}{S_t}(K-\frac{1}{T}\int_{0}^{t}S(t')dt').
\end{align}
Then, it has been shown \cite{rogers1995value,ingersoll1987theory} that the price of the Asian option satisfies the following nonlinear backward partial differential equation
\begin{equation}\label{AsianPDEs}
\left\{\begin{split}
    &\frac{\partial U}{\partial t}+\frac{\sigma^2}{2}x^2\frac{\partial U^2}{\partial x^2}- (\frac{1}{T}+rx)\frac{\partial U}{\partial x}=0,\\
     \\
     &U(T,x)=(-x)^{+},
\end{split}\right.
\end{equation}
where the partial differential equation \ref{AsianPDEs} is one dimensional PDE in space. Eq. (\ref{AsianPDEs}) is a nonlinear partial differential equation of convection-diffusion type with a convection term that is a function of volatility and spatial variable $x$. Thus, this differential equation belongs to the family of convection dominant problems where the coefficient of the convection term (second-order derivative) can be a very small number in this model. As we mentioned earlier in this section, this small coefficient could imply an oscillatory behavior such that it can lead to numerical instability for the numerical scheme \cite{douglas1982numerical}. On the other hand, the efficiency and robustness of the DPG method for the convection-diffusion problem have been successfully shown for this family of problems (\cite{ellis2016robust}, \cite{chan2014robust}, \cite{chan2013dpg} and the references therein). Demkowicz et.al. in \cite{demkowicz2013robust}, as an illustration, thoroughly analyzed the DPG method for the convection-dominated problems. They show that it benefits from a robust $L^2$ error estimate for trail variables in this set of differential equations.

Having the solution of Equ. (\ref{AsianPDEs}), the value of an Asian option with strike price $K$ and initial stock value $S_0$ can be computed as $S_0U(0,\frac{K}{S_0})$. After using a truncated computational domain $x\in [-2,2]$ for the Equ. (\ref{AsianPDEs}) and change of variable $\tau = T-t$ in time, the system of partial differential equation (\ref{AsianPDEs}) will build into the following form,
\begin{equation}\label{AsianPDEs2}
\left\{\begin{split}
    &\frac{\partial U}{\partial \tau}-\frac{\sigma^2}{2}x^2\frac{\partial U^2}{\partial x^2}+ (\frac{1}{T}+rx)\frac{\partial U}{\partial x}=0,&\forall x\in [-2,2],\quad \forall \tau \in [0,T],\\
     \\
     &U(0,x)=(-x)^{+},\\
     \\
     &U(\tau,2) = 0,\\
     \\
    &\frac{\partial U^2}{\partial x^2}(\tau,-2)= 0.\\
     \end{split}\right.
\end{equation}
So, the option value will be $S_0 U(T,\frac{K}{S_0})$. Using our convention for the DPG method, we can write the weak form for the Equ. (\ref{AsianPDEs2}) as following  
\begin{align}\label{AsianDPG}
         \frac{d}{d\tau}(u, {\bf v})+ b^{\tau}(u,{\bf v})= 0,  \\
\end{align}
where the bilinear form for primal formulation defies as 
\begin{align}\label{primalAsian}
%\left\{
\begin{split}
    b^{\tau}_{\text{primal}}({\bf u},v) = 
    (\frac{\sigma^2}{2}x^2\frac{\partial}{\partial x} u,\frac{\partial}{\partial x} v)_{\Omega} +\big((\frac{1}{T}+(r+2)x\big)\frac{\partial}{\partial x} u,v)_{\Omega}
    -\langle \hat{q}, v\rangle_{\partial \Omega},\\
 \end{split}
%\right.
\end{align}
with a set of trial and flux variables ${\bf u}=(u,\hat{q})\in H(\Omega)\times H^{\frac{1}{2}}(\Omega)$, test variable ${\bf v} = v \in L^2(\Omega)$. Moreover, considering trail variables ${\bf u} = (u,\vartheta)\in L_2(\Omega)\times L_2(\Omega)$, and flux variables ${\bf \hat{u}} = (\hat{u},\hat{{\vartheta}})\in  H^{1/2}(\Omega)\times  H^{1/2}(\Omega) $ the bilinear form (\ref{AsianDPG}) for the ultraweak formulation reads 
\begin{align}\label{ultraWkAsian}
%\left\{
\begin{split}
&b^{\tau}_{\text{ultraweak}}(({\bf u},\hat{{\bf u}}),{\bf v}) = 
    b^{\tau}_{\text{ultraweak}}(((u,{\vartheta}), (\hat{u},\hat{{\vartheta}})), (v,\omega)),
    \\
    \\
   &\qquad=-({\bf \vartheta},\frac{\sigma^2}{2}x^2 \frac{\partial}{\partial x} v)_{\Omega} +
    ({\bf \vartheta},(\frac{1}{T}+(r-\sigma^2) x  v)_{\Omega}
    -(u,\frac{\partial \omega}{\partial x})_{\Omega}
    -({\bf \vartheta}, \omega)_{\Omega}\\
    \\
   &\qquad+\langle \hat{u}, \omega\rangle_{\partial \Omega}+
    \langle \hat{{\vartheta}}, v \rangle_{\partial \Omega}.
 \end{split}
\end{align}
Now, using backward Euler approximation for time derivative and appropriate discrete test and trial space for DPG explained in the section \ref{AmericanOptionsec}, the discrete DPG formulation for the Asian option pricing problem reads
\begin{align}\label{AsianWeak}
        (u^{n+1}-u^n, {\bf v})+ \Delta \tau b^{\tau}_n(u^n,{\bf v})= 0, \qquad \forall {\bf v} \in V_h. 
\end{align}
We propose the following graph norm for ultraweak formulation and energy norm for primal DPG formulation to solve the valuing Asian option problem formulated by Equ. \ref{AsianWeak} 
\begin{align}\label{testNormAsian}
\begin{split}
         &\text{Primal}: \|v\|^2_{V} = \frac{1}{\Delta t} \|v\|^2+\frac{1}{(\Delta t)^2} \|{\sigma^2}\frac{\partial }{\partial x}v\|^2,  \\
         \\
         &\text{Ultraweak}: \|{\bf v}\|^2_{V} = \|(v,\omega)\|^2_{V}\\
         \\
         &\qquad \qquad =\frac{1}{(\Delta t)^2} \|{\sigma^2}\frac{\partial }{\partial x} v-(\frac{1}{T}+(r-\sigma^2)v-\omega\|^2+ \frac{1}{\Delta t}\|\frac{\partial }{\partial x}\omega\|^2.  \\ 
    \end{split}
\end{align}
Therefore, one can obtain the corresponding discrete operators 
\begin{align}
    B_{ij}&=b^{\tau}(u_i,{\bf v}_j),  &G_{ij}&=({\bf v}_i,{\bf v}_j)_v , &l_i &= (u_i,v).
\end{align}
However, it is worth mentioning that the above rectangle matrix $B$ is a function of the spatial variable, and the induced inner product $(\cdot,\cdot)_v$ is formed by the associated norms (\ref{testNormAsian}) defined in the procedure of the DPG formulation.
Thus, discrete DPG formulation of the equation (\ref{AsianPDEs2}) $\forall n \in \{1,\cdots, N_{\tau}\},$ yields
\begin{equation}\label{discreteAsian1}
    \left\{\begin{split}
        &B^TG^{-1}l(u^{n+1}-u^n)+ \Delta \tau B^TG^{-1}B u^{n}=0, \\
        \\
        &\hat{u}^{0} = (-x)^{+}, \qquad \forall x\in[-2,2] \\
        \\
        &\hat{u}^n|_{x=2} = 0, \qquad \forall n \in {1,\cdots, N_{\tau}},    
    \end{split}\right.
\end{equation}
Thus, we can define near the optimal discrete operators  $B^{\text{n-op}}=B^T G^{-1}B$,$l^{\text{n-op}}=B^T G^{-1}l$ discrete DPG for the equ (\ref{AsianPDEs2}) for all  $\forall n \in \{1,\cdots, N_{\tau}\},$ finds 
\begin{equation}\label{discreteAsian2}
    \left\{\begin{split}
        &l^{\text{n-op}}(u^{n+1}-u^n)+ \Delta \tau B^{\text{n-op}} u^{n}=0, \\
        \\
        &\hat{u}^{0} = (-x)^{+}, \quad \forall x\in[-2,2], \\
        \\
        &\hat{u}^n|_{x=2}  = 0, \qquad \forall n \in \{1,\cdots, N_{\tau}\}.
    \end{split}\right.
\end{equation}
The system of Equ. (\ref{discreteAsian2}) can be solved by an iterative method or linear solver. In the next section, we examine the efficiency of the proposed DPG method.

\subsection{Numerical Experiments}
As mentioned before, the set of the partial differential equations (\ref{AsianPDEs2}) is a nonlinear and convection-dominant problem, and developing a numerical scheme for this problem can be problematic due to the convection term.
In this section, we select some famous test problems from the literature to showcase the efficiency and accuracy of the proposed numerical scheme (\ref{discreteAsian2}). In this example, all the results are generated by the first-order DPG method, and corresponding to the enriched test spaces ($\Delta p = 2$). we used $N_s=100$ number of spatial elements, and the $N_t=100$ time step for all the experiments in this section.
%%%%%%%%%%%%%%%%%%%%%%%%%%%%%%%%%%%%%%%%%%%%%%%%%%%%%%%
%-----------------Value of Asian option----------------
%%%%%%%%%%%%%%%%%%%%%%%%%%%%%%%%%%%%%%%%%%%%%%%%%%%%%%%
\begin{figure}[!ht]
\includegraphics[width=\linewidth]{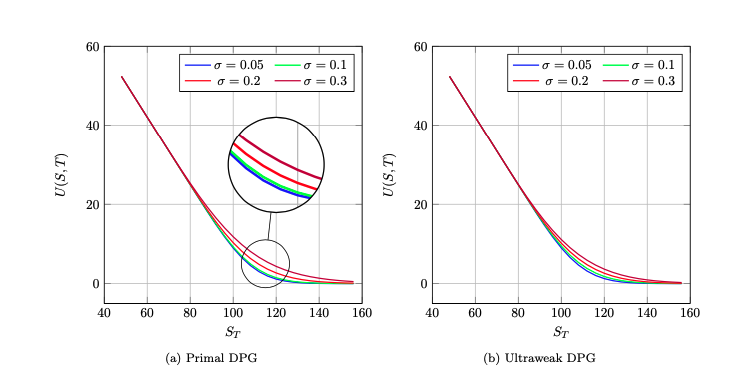}
\caption{Value of an Asian option with DPG method for $r = 0.015$, and different volatility }
\label{DPGAsianOption}
\end{figure}
Fig. (\ref{DPGAsianOption}) displays the value of the Asian option with two ultraweak and primal DPG formulations for different values of $\sigma = 0.05,0.1,0.2, 0.3$.
As it can be seen the value of the Asian option is smooth and stable even for a small value of $\sigma = 0.05$ which leads to the convection-dominated case for the system of \ref{discreteAsian2}.
%----------------------------------------------------------------------------
%%%%%%%%%-------------------------Table Asian option -----------------------
%----------------------------------------------------------------------------
\begin{table}[htbp]
\footnotesize
  \caption{{Asian call option with $r=0.09$, $T=1$, $S_0=100$ }\label{AsianTable1}}
   \centering
 \scalebox{0.85}{ 
  \begin{tabular}{c|l|l|c|c|c} \hline
   $\sigma$ & Reference &Method& $K=95$&$K=100$&$K=105$\\ \hline
\multirow{7}{*}{0.05}& Zhang \cite{zhang2001semi} &   & 8.8088392 & 4.3082350 & 0.9583841\\ 
            & Zhang-AA2 \cite{zhang2003pricing}   &   &8.80884 &4.30823 &0.95838\\ 
            & Zhang-AA3 \cite{zhang2003pricing}  &  & 8.80884 & 4.30823 & 0.95838\\ 
            &                         & Ultraweak DPG&  8.8085332 & 4.3081967 & 0.958371\\ 
            &                         & Primal DPG   &  8.8088363 & 4.3082291 & 0.9583836\\
\hline
\multirow{7}{*}{0.10}& Zhang\cite{zhang2001semi}   &   & 8.9118509 &4.9151167 & 2.0700634\\ 
            &         Zhang-AA2\cite{zhang2003pricing}&  & 8.91171& 4.91514 & 2.07006\\ 
            &      Zhang-AA3\cite{zhang2003pricing}  &    & 8.91184 & 4.915126 & 2.07013\\  
            &                              & Ultraweak DPG& 8.910986 &4.915116769 & 2.0700633\\ 
            &                              & Primal DPG   & 8.9118498 &4.9151265 & 2.0700634\\
\hline
\multirow{7}{*}{0.20}&    Zhang \cite{zhang2001semi}&  & 9.9956567 & 6.7773481 & 4.2965626\\ 
            &      Zhang-AA2\cite{zhang2003pricing} &&9.99597 &6.77758 & 2.745\\ 
            &    Zhang-AA3\cite{zhang2003pricing} &  & 9.99569 & 6.77738 & 4.29649\\  
            &                          & Ultraweak DPG& 9.99565668 & 6.7773481 & 4.2965626\\  
            &                          & Primal DPG  & 9.9956567 & 6.7773429 & 4.2965619\\  
\hline
\multirow{7}{*}{0.30}&Zhang\cite{zhang2001semi} &   &11.6558858&8.8287588 &6.5177905\\ 
            &Zhang-AA2 \cite{zhang2003pricing}  &   &11.65747&8.82942& 6.51763\\ 
            &Zhang-AA3 \cite{zhang2003pricing} &  & 11.65618&8.82900 &6.51802\\  
            &                          & Ultraweak DPG&11.6558853&8.8287498 &6.51779047\\
            &                          & Primal DPG   &11.6558857&8.8287580 &6.51779054\\ 
            
\hline
  \end{tabular}
}
\end{table}
%%%%%%%%%%%%%%%End Table Asian option %%%%%%%%%%%%%%%%
we prepared the table \ref{AsianTable1} to compare the result of DPG methodology for pricing an Asian option with interest rate $r=0.09$, $T=1$, $S_0=100$, different strike price $K=95,100,105$, and different volatility with the result report in \cite{zhang2001semi}, \cite{zhang2003pricing}. Considering the result from \cite{zhang2001semi} as a benchmark with the PDE method, one can see that the obtained results from DPG ultraweak and primal method are so close (less than $0.001\%$ deviation).
%%%%%%%%%%%%%%%%%%%%%%%%%%%%%%%%%%%%%%%%%%%%%%%%%%%
%Rogers and Shi (1995)\cite{rogers1995value}
% Foufas and Larson (2008) \cite{foufas2008valuing}
\begin{table}[htbp]
\footnotesize
  \caption{{Asian call option with $r=0.15$, $T=1$, $S_0=100$ }\label{AsianTable2}}
   \centering
 \scalebox{0.85}{ 
  \begin{tabular}{c|l|l|c|c|c} \hline
   $\sigma$ & Reference &Method& $K=95$&$K=100$&$K=105$\\ \hline
\multirow{7}{*}{0.05}    & vecer \cite{vecer2001new} & Monte Carlo  & 11.094 & 6.795 & 2.745\\ 
            &      & FDM          & 11.094 & 6.795 & 2.745\\ 
            & Rogers and Shi \cite{rogers1995value}    & Lower Bound  & 11.094 & 6.795 & 2.745\\ 
            & Foufas and Larson \cite{foufas2008valuing} & FEM          & 11.112 & 6.810 & 2.754\\ 
            & Kim et al. \cite{kim2007meshfree}         & MPCM         & 11.093 & 6.79 & 2.78\\ 
            &                          & Ultraweak DPG& 11.09398 & 6.79512 & 2.74481\\ 
            &                          & Primal DPG   & 11.09401 & 6.7948 & 2.74499\\ 
\hline
\multirow{7}{*}{0.10}& vecer\cite{vecer2001new}      & Monte Carlo  & 15.399 & 7.028 & 1.418\\ 
            &        & FDM          & 15.399 & 7.029 & 1.415\\ 
            & Rogers and Shi \cite{rogers1995value}   & Lower Bound  & 15.399 & 7.028 & 1.413\\ 
            & Foufas and Larson \cite{foufas2008valuing} & FEM          & 15.416 & 7.042 & 1.422\\ 
            & Kim et al. \cite{kim2007meshfree}         & MPCM         & 15.398 & 7.028 & 1.448\\ 
            &                          & Ultraweak DPG& 15.3984 & 7.0277 & 1.41769\\
            &                          & Primal DPG   & 15.39899 & 7.02812 & 1.418001\\
\hline
\multirow{7}{*}{0.20}    & vecer\cite{vecer2001new}  & Monte Carlo  & 15.642 & 8.409 & 3.556\\ 
            &                          & FDM          & 15.643 & 8.412 & 3.560\\ 
            & Rogers and Shi\cite{rogers1995value}    & Lower Bound  & 15.641 & 8.408 & 3.554\\ 
            & Foufas and Larson \cite{foufas2008valuing} & FEM          & 15.659 & 8.427 & 3.570\\ 
            & Kim et al. \cite{kim2007meshfree}         & MPCM         & 15.66437 & 8.421 & 3.573\\ 
            &                          & Ultraweak DPG& 15.64218 & 8.4091 & 3.5559\\ 
            &                          & Primal DPG   & 15.641865 & 8.4102 & 3.5584\\
\hline
\multirow{7}{*}{0.30}& vecer\cite{vecer2001new}     & Monte Carlo  & 16.516 & 10.210 & 5.731\\
            &                          & FDM          & 16.516 & 10.215 & 5.736\\ 
            & Rogers and Shi\cite{rogers1995value}    & Lower Bound & 16.512 & 10.208 & 5.728\\ 
            & Foufas and Larson\cite{foufas2008valuing} & FEM          & 16.553 & 10.231 & 5.750\\ 
            & Kim et al. \cite{kim2007meshfree}        & MPCM         & 16.5179 & 10.2194 & 5.742\\ 
            &                          & Ultraweak DPG& 16.51615 & 10.21045 & 5.73074\\ 
            &                          & Primal DPG   & 16.51617 & 10.20964 & 5.730865\\ 
          \hline
  \end{tabular}
}
\end{table}
%%%%%%%%%%%%%%%End Table Asian option %%%%%%%%%%%%%%%%
To compare the accuracy and stability of the proposed method with the broader method in the literature, table. (\ref{AsianTable2}) is produced. In this test, the results from the Monte Carlo method are exploited as an exact solution. we compute the value of an Asian option for different strike prices $K=95,100,105$, the interest rate of $r = 0.15$, time to maturity $T=1$, initial asset value $S_0=100$, with different volatility $\sigma = 0.05,0.1,0.2,0.3$. The result from the DPG methods is a maximum $0.001\%$ deviation from the benchmark. 
%----------------------End of asian option ----------------------------------------
%%---------------------------------------------------------------------------------
%-------------------------------------Barrier Option--------------------------------
\subsection{Barrier Options}
A double knock-out Barrie option is a financial contract that gives a payoff $h(S)$ at maturity $T$, as far as the price of the underlying asset stays in the predetermined barriers $[S_L(t), S_U(t)]$, otherwise, if the spot price is hit barriers, the option gets knocked out. Although the barriers are checked continuously in time, it is more feasible to check the barriers discretely in the real-world application \cite{shreve2004stochastic}.

It is well-known that the closed-form analytical solution for the discrete double barrier option is not known, so devising accurate and efficient numerical methods for valuing this type of option is essential. Thus, over the past years, researchers try to develop semi-analytical and numerical schemes for approximating the price of Barrier options. Here, we briefly address some of them. Kunitomo et.al \cite{kunitomo1992pricing} used sequential analysis to find the solution as a series, analytical approach by contour integration is used by Pelsser \cite{pelsser2000pricing} to price the barrier options. The binomial method is used by Cheuk et. al in \cite{cheuk1996complex}, and  the Monte Carlo method as a probability-based method is devised in \cite{ndogmo2007high} to price this exotic option
PDE method such as the finite difference method by Zevan et. al in \cite{zvan2000pde}, a finite element in \cite{golbabai2014highly} by Golbabai et.al, and quadrature method in \cite{milev2010numerical}  is developed for the pricing discrete barrier options.

We begin by stating the model of the problem which is inspired by the work   \cite{milev2010numerical}, and \cite{tse2001pricing}. Assume that dynamic of the underlying asset $\{S_t\}_{t\in [0,T]}$ is following the stochastic differential equation in  (\ref{geoBrownMo}), with the standard winner process $W_t$, interest rate $r$, volatility of $\sigma$, and fixed initial asset price $S_0$. Defining the Brownian motion $Z_t$, with instantaneous drift value $\hat{r} = r-(\sigma^2/2)$, and standard deviation $\sigma$, the price process will follow $S_t = S_0 e^{Z_t}$. Moreover, we define the discrete counterpart process $\Tilde{X}_n = S_0 e^{\Theta_n}$, for $n=1,2,\cdots,N$, and $\Theta_n=\theta_1+\theta_2+\cdots+\theta_n$, $\Theta_0 = 0$. Random variables $\theta_i$ are independent normally distributed random variables i.e. $N(r-\sigma^2/2, \Delta t \sigma)$ with $\Delta t = \frac{T}{N}$ for $N$  predetermined monitoring instants. 

Consider the discrete monitoring dates of $t_1=0\leq t_2\leq \cdots \leq t_N=T$ with the constant upper and lower barriers of $S_U$, and $S_L$ respectively. Besides, we assume that barriers are not active on the first, and last dates of our time interval. The price of a discrete double barrier option can be computed by discount of expected payoff at expiration time $T$ to the present time $t$ as follows.
\begin{align*}
   e^{-rT} E [h(S_T)| \chi_{B_1}\chi_{B_2}\cdots \chi_{B_n}],
\end{align*}
where the indicator functions of $\chi_{\cdot}$ is evaluating on sub set of $B_i=\{S_i\in(S_L,S_U)\}$. 

Denoting $U(t, S)$ the value of a discrete double barrier option with the date of maturity of $T$, strike price $K$,(for simplicity), this value will satisfy in the following system of $N$ partial differential equations
\begin{equation}\label{barrierPDE1}
\left\{
\begin{split}
          &\frac{\partial U(t,S)}{\partial t}+\frac{\sigma^2}{2}S^2\frac{\partial U(t,S)^2}{\partial S^2}
          +r S\frac{\partial U(t,S)}{\partial S}-rU(t,S)=0,
          \hspace{0.0cm} \quad \forall t\in [t_i,t_i+1], \quad \forall i=1,2,\cdots,N \\
          \\
         &U(t,0 )=0,\\
         \\
         &\lim_{S\to +\infty}U(t,s) = h(S)\\
         \\
         &U(S,t_i) = h_i(S)\qquad \forall i=1,2,\cdots,N\\
         \\
         &U(S,T) = h_T(S), 
\end{split}
\right.
\end{equation}
where boundary conditions $h_i(S)$, and $h_T(S)$ are also defined as
\begin{equation}\label{barrierBC1}
h_i(S) = \left\{
         \begin{split}
         &\lim_{t\to t^{+}} U(S,t),\quad \text{if}\quad  S_L\leq S\leq S_U,       \hspace{0.5cm}\quad \forall t\in [t_i,t_i+1], \quad \forall i=1,2,\cdots,N\\
         \\
            &  0, \hspace{2.5cm} \text{if}\quad S=\mathbb{R}^{+}\backslash [S_L,S_U], 
         \end{split}
         \right.
\end{equation}
,and 
\begin{equation}\label{barrierBC2}
h_T(S) = \left\{
         \begin{split}
          &(S-K)^{+},\quad \text{if}\quad  S_L\leq S\leq S_U,       \hspace{0.5cm}\quad \forall t\in [t_i,t_i+1], \quad \forall i=1,2,\cdots,N,\\
         \\
           &   0,  \hspace{2.5cm}  \text{if} \quad S=\mathbb{R}^{+}\backslash [S_L,S_U]. 
         \end{split}
         \right.
\end{equation}
As we can observe, the set of partial differential equations (\ref{barrierPDE1}) is a system of consecutive partial differential equations where on each time interval $[t_i,t_{i+1}]$ has the final boundary conditions of (\ref{barrierBC1}), and the final PDE has the boundary condition (\ref{barrierBC2}). Besides, the system of PDEs presented in (\ref{barrierPDE1}) with the aforementioned boundary condition is a non-smooth and nonlinear partial differential equation associated with the functions (\ref{barrierBC1}), (\ref{barrierBC2}), therefore, designing an accurate and stable numerical scheme is tricky here. 

We use the change of variable in space and time similar to the change of variables for vanilla options in section (\ref{BlkSholes}) to obtain the following piecewise constant coefficient partial differential equations.
\begin{equation}\label{barrierPDE}
    \left\{
    \begin{split}
          &\frac{\partial u(\tau,x)}{\partial \tau}+\frac{\sigma^2}{2}\frac{\partial u(\tau,x)^2}{\partial x^2}
          +(r+\frac{\sigma^2}{2}) \frac{\partial u(\tau,x)}{\partial x}-ru(\tau,x)=0,
          \hspace{0.0cm}\quad \forall \tau\in [\tau_i,\tau_{i+1}], \quad \forall i=1,2,\cdots,N \\
          \\
         &u(\tau,0 )=0,\\
         \\
         &\lim_{x\to +\infty}u(\tau,x) = h(x),\\
         \\
         &u(x,\tau_i) = h_i(x),\qquad \forall i=1,2,\cdots,N,\\
         \\
         &u(x,T) = h_T(x). 
          \end{split}
    \right.
\end{equation}

Now if we concentrate on one of the equations as a generic differential equation on the interval $[\tau_j,\tau_{j+1}]$, where $j\in \{1,2,\cdots,N\}$, we propose the following weak formulation for DPG formulation 
\begin{align}\label{BarrirDPG}
            \frac{d}{d\tau}({\bf u}, {\bf v})+ b^{\tau}({\bf u},{\bf v})=0,
         \qquad \forall \tau\in [\tau_j,\tau_j+1],
\end{align}
where the bilinear form is similar to the primal and ultraweak formulation defined in Eq. (\ref{primalVI}), and (\ref{ultraWkVI}) on this sub-interval. However, the boundary conditions introduced (\ref{barrierPDE}) are performing on the interval $[\tau_j,\tau_{j+1}]$ as a sub-interval of the computational domain. Utilizing a generic partition $\tau_{j}=\tau_{j1},\tau_{j2},\cdots, \tau_{jN_j} =\tau_{j+1},$ for each interval, and using backward Euler scheme for time derivative, the approximate of equation (\ref{BarrirDPG}) in the finite dimension space, the discrete DPG for each sub-partial differential equations reads
\begin{align}\label{discrtBarier}
            (u^{n+1}-u^n, {\bf v})+ \Delta \tau_i b^{\tau}_n(u^n,{\bf v})= 0, \qquad \forall {\bf v} \in V_h., \quad \forall n\in \{1,2, \cdots, N_{i}\},
\end{align}
where the time steps on the domain of each sub-problem defined as $\Delta \tau_i = \frac{\tau_{i+1}-\tau_i}{h_i},\quad h_i = \frac{\tau_{i+1}-\tau_i}{N_i}$. Indeed, on each problem (\ref{discrtBarier}) we need to solve a nonlinear non-smooth discrete system of equations (see the psudo code \ref{algBarrier}).  Defining the graph and energy norm defined in (\ref{testNormVanial}) for each sub-domain $[\tau_j,\tau_{j+1}]$, and denoting the discrete operators of  $B$, $G$, and $l$ accordingly as following
\begin{align}
    B_{ij}&=b^{\tau}(u_i,{\bf v}_j), & G_{ij}&=({\bf v}_i,{\bf v}_j)_v , & l_i &= (u_i,v), & U_h&=[\tau_j,\tau_{j+1}],
\end{align}
One can find the discrete nonlinear generic problems on each sub-domain
\begin{equation}\label{discreBarreri2}
    \left\{\begin{split}
         &l^{\text{n-op}}(u^{n+1}-u^n)+ \Delta \tau B^{\text{n-op}} u^{n}=0,\qquad  \forall n \in \{1,\cdots, N_{\tau_i}\}, \\
         \\
         &{u}^{n}|_{x=0}=0, \qquad \forall n \in \{1,\cdots, N_{\tau_i}\},\\
         \\
         &\lim_{x\to +\infty}{u}^{n+1}(x) = h(x),\\
         \\
         &{u}^{{N_{\tau_i}}}|_x = h_i(x).\\
            \end{split}\right.
\end{equation}
where the near optimal DPG operators are denied as $B^{\text{n-op}}=B^T G^{-1}B$, $l^{\text{n-op}}=B^T G^{-1}l$. The non-smooth system of equations of (\ref{discreBarreri2}) can be solved by a projected iterative solver such as Gradient descent for different consecutive intervals till the time of maturity \cite{beck2014introduction}.
%%%%%%%%%%%%%%%%%%%%%%%%%%%% Algrithm %%%%%%%%%%%%%%%%%%%%%%%%%%%%
\begin{algorithm}
\caption{Numerical algorithm for the double barrier option}\label{algBarrier}
\begin{algorithmic}
\Require $S_0\in[S_L,S_U]$
\State $u^N|_S \gets h(S)$
\For{$\tau_i\in[t_1,t_N]$}
\For{$\tau_{i_j}\in[\tau_i,\tau_{i+1}]$}
\If{$S$ is in $[S_L,S_U]$}
\State $u^n|_{S=0} = 0$,
\State $u^{N_{\tau_i}}|_{S\to\infty} = h(S)$,
\State $u^{N_{\tau_i}}|_S = h_i(S)$,\\
\State Solve the sub-partial differential equation \ref{discreBarreri2}.\\
\ElsIf{$S$ is out of $[S_L,S_U]$}\\ 
\State The option will be knocked out!\\
\EndIf
\EndFor
\EndFor
%\State $y \gets 1$
%\State $X \gets x$
\end{algorithmic}
\end{algorithm}

%%%%%%%%%%%%%%%%%%%%%%%%%%%%End Algrithm %%%%%%%%%%%%%%%%%%%%%%%%%%%%
\subsection{Numerical Experiments}
Here we solve the standard test problem solved in \cite{kim2014option} problem. we use the DPG method to price a barrier option with volatility $\sigma=0.2$, interest rate $r = 0.1$, strike price $k=100$, and upper and lower boundary of $S_L=95$, and $S_U = 125$ respectively. It is known that a trading year includes $250$ a working day, and a working week has five days. In this example, we report the numerical estimate for daily and weekly monitoring. in another words, if we take $T=1$ (half year $T=0.5$) for one trading year , then time increments of $\Delta t = 0.004$ (half year $\Delta t = 0.002$) corresponds with daily check and $\Delta = 0.02$ (half year $\Delta t = 0.01$) corresponds to weekly check. Using the first-order DPG method with $N_s=100$ spatial element, $N_t = 100$ stepping time, and enriched test space with $\Delta p=2$, the desired results will accomplish.     

\begin{figure}[!htb]
\includegraphics[width=\linewidth]{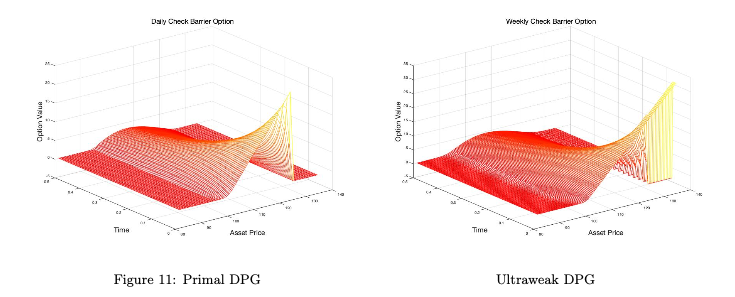}
\caption{Surface of the price of barrier option, $\sigma =0.2$, r =0.1, K=100, $[S_L,S_p]=[95,125]$ via DPG method. }\label{UltraPrimaPriceBarrier}.
\end{figure}

\begin{figure}[!ht]
\includegraphics[width=\linewidth]{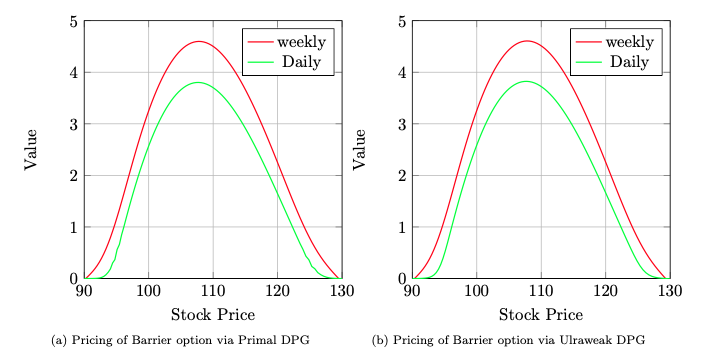}
\caption{Value of Barrier option with DPG method}
\label{DPGBbarrierValue}
\end{figure}

Fig. \ref{UltraPrimaPriceBarrier} depicts the surface of the price of the barrier option with the two primal and ultraweak formulations, as we expect this option is cheaper than the European option due to the convenience that brings for the trader. Moreover, in spite of the non-smooth boundary condition the surface of the price is smooth and stable.  
%-------------------------- price surface primal and ultraweak----------------
%-------------------------- End price surface primal and ultraweak----------------
We prepared fig. \ref{DPGBbarrierValue} to show the price of the barrier option with the aforementioned market parameters. The primal and ultraweak formulation is implemented to find the value of the option by checking both weekly and Daily for the barriers. One difficulty in pricing barrier options is that the value of the option can be oscillatory near the barriers of $S_L$, and $S_U$, whereas the illustrations show the stable and smooth behavior of the price for the value of stock price close to the boundaries.
%----------------------------price of Barrier option------------
%---------------------------------End price Barrier option---------------------------
Table (\ref{BarrierTable}) compares the accuracy of the DPG method with the path integral method \cite{milev2010numerical}, and MPCM method \cite{kim2007meshfree}. in this experiment, we see the value of the option for daily and weekly monitoring when the price of the underlying price is $S=95$, $S=95.0001$, $S=124.9999$, and $S=125$ when barriers are $[S_L, S_U] = [95,125]$. As we mentioned before the numerical scheme can have unstable behavior close to barriers, and in this example, we try to catch the accuracy of the method when the stock price is in a very close neighborhood of barriers. As can be seen, the DPG method is accurate and very close to the recorded value in \cite{milev2010numerical}, and \cite{kim2007meshfree}.
%---------------------------------table Barrier option---------------------------
{\small 
\begin{table}[H]
{\footnotesize
  \caption{{Double Barrier option with $\sigma=0.2$, $r=0.1$, $T=0.5$, $K=100$, $L=95$,$U=125$ }\label{BarrierTable}}
\begin{center}
  \begin{tabular}{c|l|l|c|c} \hline
   $S$ & Reference &Method& Weekly checking & Daily checking\\ \hline
\multirow{4}{*}{95}
& Milev \cite{milev2010numerical} & Path integral& 11.094 & 6.795 \\ 
& Kim \cite{kim2007meshfree}    & MPCM         & 11.094 & 6.795 \\ 
                   &        & Primal DPG   & 11.094 & 6.795 \\ 
                   &       & Ultraweak DPG& 11.094 & 6.795 \\ 
\hline 
\multirow{5}{*}{95.0001}
& Milev \cite{milev2010numerical} & Path integral& 11.094 & 6.795 \\ 
              & Monte Carlo & $10^7$ paths & 11.094 & 6.795 \\
&Kim \cite{kim2007meshfree}     & MPCM         & 11.094 & 6.795 \\ 
                   &        & Primal DPG   & 11.094 & 6.795 \\ 
                   &        & Ultraweak DPG& 11.094 & 6.795 \\ 
\hline 
\multirow{5}{*}{124.9999}
&Milev \cite{milev2010numerical}  & Path integral& 11.094 & 6.795 \\ 
              & Monte Carlo & $10^7$ paths & 11.094 & 6.795 \\
&Kim \cite{kim2007meshfree}     & MPCM         & 11.094 & 6.795 \\ 
                   &        & Primal DPG   & 11.094 & 6.795 \\ 
                   &        & Ultraweak DPG& 11.094 & 6.795 \\ 
\hline 
 \multirow{4}{*}{125}
 & Milev\cite{milev2010numerical} & Path integral& 11.094 & 6.795 \\ 
&Kim \cite{kim2007meshfree}      & MPCM         & 11.094 & 6.795 \\ 
                   &         & Primal DPG   & 11.094 & 6.795 \\ 
                   &         & Ultraweak DPG& 11.094 & 6.795 \\ 
\hline            

  \end{tabular}
\end{center}
}
\end{table}
}
%-----------------------------END-table Barrier option---------------------------
%%%%%%%%%%%%%%%%%%%%%%%%%%%%%%%%%%%%%%%%%%%%%%%%%%%%%%%%%%%%%%%%%%%%%%%%%%%%
%---------------------------- THE END of Barrier option --------------------
%%%%%%%%%%%%%%%%%%%%%%%%%%%%%%%%%%%%%%%%%%%%%%%%%%%%%%%%%%%%%%%%%%%%%%%%%%%%
%-------------------------------Section Analysis-------------------------------
\section{Option Pricing problem with Stable Method }

\begin{thm}[Main Theorem]\label{thm:Mainthm}
option pricing with DPG has a unique solution!!!
\end{thm}
%%%%%%%%%%%%%%%%%%%%%%%%%%%%%%%%%%%%%%%%%%%%%%%%%%%%%%%%%%%%%%%%%%%%%%%%%%%%
%---------------------------- Sensitivity Analysis ------------------------
%---------------------------------------------------------------------------
%-------------------------------Section Greeks-------------------------------
\section{Sensitivity Analysis with Greeks}\label{greeks}
In this section, we use the DPG methodology to calculate the sensitivity of option pricing under the Black-Scholes model. Sensitivity of the option with respect to model parameters, Greeks, explains the reaction of the option value to the fluctuation of the market environment. Greeks are compasses in the trader's hand to find the correct direction in the hope of hedging their portfolio by buffering against market changes.
Thus, the efficiency and accuracy of the numerical scheme are of paramount importance to trace the option price changes when the state of the market changes. Let $u(s,t)$ be the solution of Black-Scholes partial differential (\ref{BlS1}) with the appropriate boundary condition pertaining to that specific option, and $\alpha$ is the desired parameter for which we want to see the changes of price, then $\frac{\partial u(x,t)}{\partial \alpha}$ which for simplicity it will be denoted by $u_{\alpha}(x,t)$ is the sensitivity. This sensitivity can be found with the direct method or dual method (the avid readers can see \cite{damircheli2019solution}). Taking the derivative with respect to the parameter $\alpha$ from Black-Scholes, one can find a system of partial differential equation that seeks for $u_{\alpha}(x,t)$
\begin{align}\label{sensPDE1}
    \frac{\partial u_{\alpha}}{\partial t}+
        \frac{\partial}{\partial \alpha}(\frac{\sigma^2}{2}x^2)\frac{\partial u^2}{\partial x^2}+\frac{\sigma^2}{2}x^2\frac{\partial u_{\alpha}^2}{\partial x^2}
        +\frac{\partial}{\partial \alpha}(r x) \frac{\partial u}{\partial x}+
        r x \frac{\partial u_{\alpha}}{\partial x}
        -\frac{\partial r}{\partial \alpha}u(x,t)- r u_{\alpha} = 0,
\end{align}
Note, the $u(x,t)$ is already evaluated the value of the option in the initial state of parameter $\alpha$ (see \cite{seydel2006tools} for more detail). One can develop a DPG formulation either primal or ultraweak for solving the PDE presented in (\ref{sensPDE1}) to find the desired sensitivity of $u_{\alpha}(x,t)$ with appropriate boundary condition. In this paper, we study the first and second derivative of price with respect to the underlying asset that are named as Delta and Gamma respectively. 

\begin{figure}[!ht]
\includegraphics[width=\linewidth]{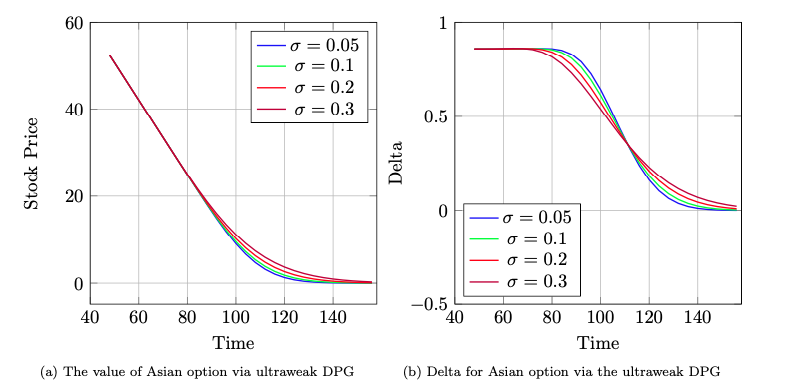}
\caption{Computing the Delta for Asian option alongside the value of the option with the ultraweak DPG method}
\label{GreekUltraAsian}
\end{figure}

\begin{figure}[H]
\includegraphics[width=\linewidth]{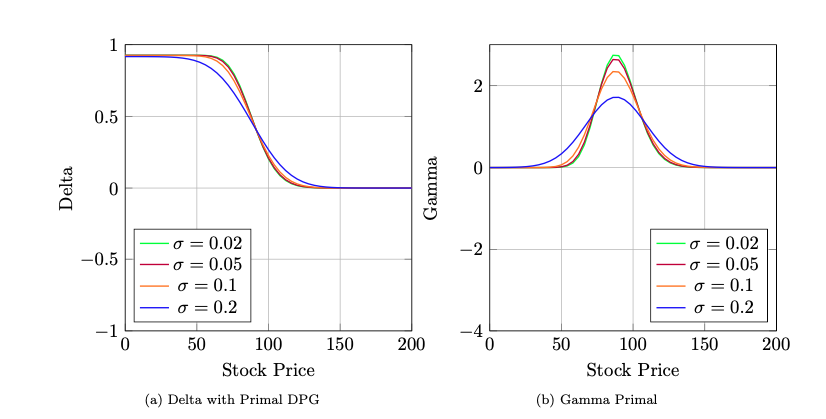}
\caption{Greeks of Asian option with primal DPG method}
\label{DirectGreekAsianPrimal}
\end{figure}

\begin{figure}[!ht]
\includegraphics[width=\linewidth]{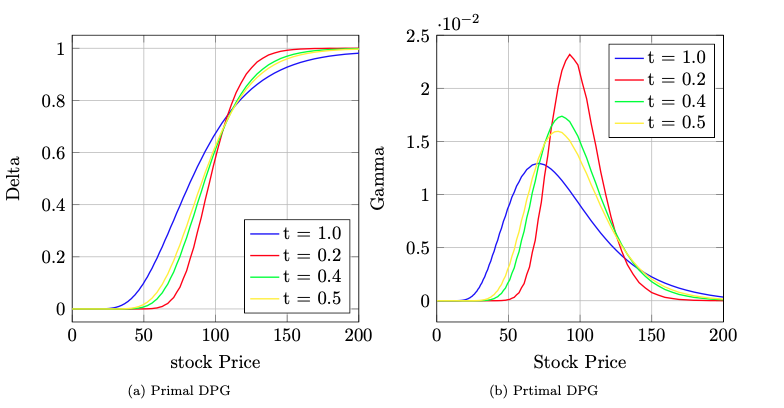}
\caption{Greek of European call option with Primal DPG parameters: r=0.05, $\sigma$=0.15,K=100}
\label{GreekEuropeanCallPrimal}
\end{figure}

\begin{figure}[H]
\includegraphics[width=\linewidth]{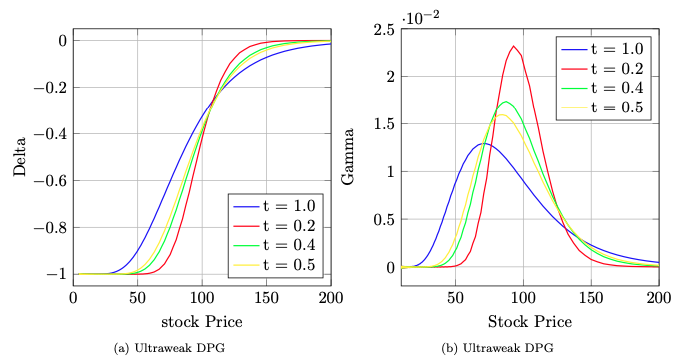}
\caption{Greek of European Put option with ultraweak DPG parameters: r=0.05, $\sigma$=0.15,K=100}
\label{GreekEuropeanPutUltrawk}
\end{figure}

\begin{figure}[H]
\includegraphics[width=\linewidth]{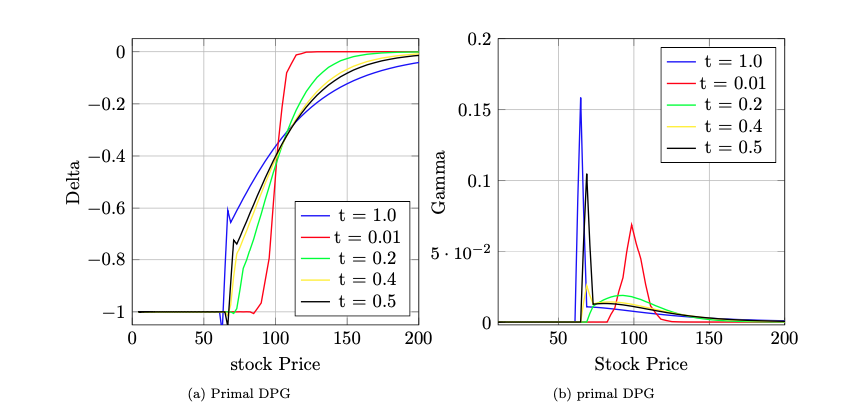}
\caption{Greek of American Put option with primal DPG parameters: r=0.05, $\sigma$=0.15,K=100}
\label{GreekAmericanPrimal}
\end{figure}

\begin{figure}[!ht]
\includegraphics[width=\linewidth]{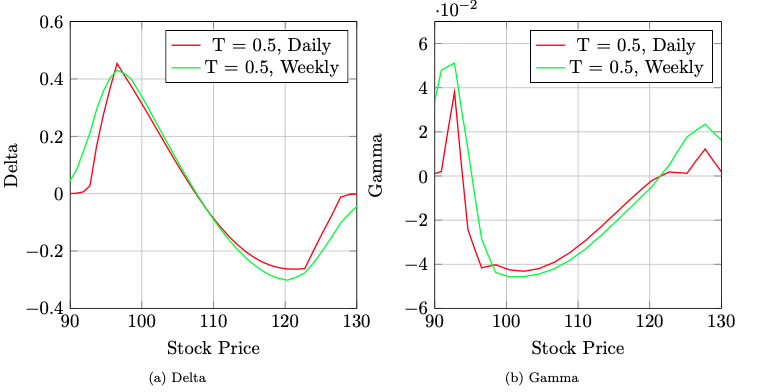}
\caption{primal Barrier for weekly and daily, p=2, Nt = 25, S0=95.0001}
\label{GreekBarPrimal}
\end{figure}

\begin{figure}[!ht]
\includegraphics[width=\linewidth]{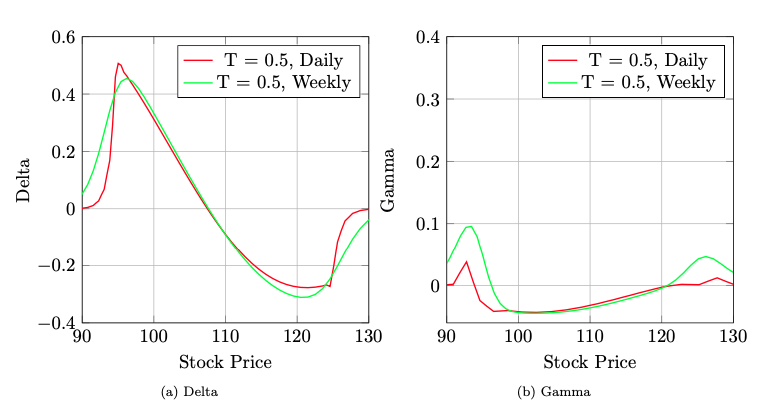}
\caption{Ultraweak Barrier for weekly and daily, p=2, Nt = 25, S0 = 100}
\label{GreekBarultraweak}
\end{figure}

To start, it is worth mentioning that in the ultraweak formulation of DPG method (for example see (\ref{ultraVanila})) inherently and implicitly we are evaluating the Delta since our primary trail variables are $(u(x,t), \sigma = \frac{\partial u}{\partial x})$. 
Fig. (\ref{GreekUltraAsian}) the numerical result of  ultraweak solution of the Asian option pricing problem as an example is prepared to show how Delta can implicitly be calculated without extra computational cost for recalculation of sensitivity.  
%%%%%%%%%%%%%%%%%%%%%%%%%%%%%%%%%%%%%%%%%%%%%%%%%%%%%%%
%-------------sol and derivative of Asian option-------
%%%%%%%%%%%%%%%%%%%%%%%%%%%%%%%%%%%%%%%%%%%%%%%%%%%%%%%
%---------------------------------End of figure Asian greek ultraweak-----------
However, one can indirectly find the Gamma and Delta of Asian option with ultraweak formulation and primal formulation Fig. (\ref{DirectGreekAsianPrimal}) using the PDE (\ref{sensPDE1}) for different volatility of the market. 
%---------------------------------------------------------------------
%%%%%%%%%%%%%% Sensitivity Primal Asian ---%%%%%%%%%%%%%%%%%%%%%%%%%%%%
%---------------------------------------------------------------------
%--------------------------------End Primal Asian Greeks-----------
%-------------------------------End of Sensitivity of Asian Option---------------
%---------------------------------------------------------------------------
It is well-known that delta is positive for call options Fig. (13a) and negative for put option Fig. (13b), whereas Gamma is always positive for both call options Fig. (14a) and put options Fig. (14b). Fig.(\ref{GreekEuropeanCallPrimal}) is prepared to illustrate Delta and Gamma of the European call option for different times to maturity, strike price $k=100$, $r = 0.05$, and $\sigma = 0.15$ with primal DPG method. The sensitivity of the European put option with the same market parameters is depicted in Fig. (\ref{GreekEuropeanPutUltrawk}) using the ultraweak DPG method.
%---------------------------EuropeanCall-----------------------------------
%---------------------------End European Call-----------------------------------
%-------------------------------------------------------------------------------
%-------------------------------------------------------------------------------
%------------------------European PUT Ultraweak---------------------------------
Admittedly, the American option is one of the most attractive options for market makers since they have the right to exercise the contract once they find the appropriate moment based on their hedging strategy. Thus, not only the monitoring Delta is important, but practitioners are curious about the rate of change in Delta itself (Gamma) for each one-basis point movement in the underlying asset. However, we can expect that the free boundary attained by the early exercise feature has a significant impact on the sensitivity of the option as well. Fig.(\ref{GreekAmericanPrimal}) shows the violation in Delta and Gamma for an American Put option based on the Primal DPG method in the different time to maturities. As we can see this chaotic behavior as the time approaches maturity increases such that at $t=0.01$ shortly after locking the option we have smooth behavior like the European option and at time $t=1.0$ we have maximum fluctuation. 
%-------------------Greek Aaerican_PUT_primal-----------------------------------
%----------------------End Greek American Put primal--------------------------
Greeks for barrier option with the double barrier $u_l = 95$, and $u_u = 125$ has shown in Fig. (\ref{GreekBarPrimal}), and Fig. (\ref{GreekBarultraweak}) using DPG method for different initial stock price $S0 = 95.0001$, $S0 = 100$. Both figures show that the sensitivity has sinusoidal behavior around the barriers when the underlying price is close to $95$, and $125$. One can see that in both cases rate of change in price and Delta are more smooth for weekly checking the barriers in comparison to daily check of the barrier which stands to reason.
%--------------------- Greeks primal Barrier Barrier primal------------

%--------------------------End Greeks primal Barrier------------------

%----------------------------Barrier Ultraweak--------------------------

%--------------End Greek ultrawek Barrier Ultraweak Barrier------------------
\section{Conclusion}
In this manuscript, a numerical scheme based on the discontinuous Petrov–Galerkin (DPG)
is proposed to deal with the option pricing problem as one of the most important branches of quantitative finance. The Black-Scholes PDE arisen from option pricing  is a special member of the family of the convection-diffusion problem which is known for being unstable in the case of having a convention-dominant term. The DPG method automatically yields a stable numerical method to estimate the solution of the very same PDE. In this investigation, we derived detailed DPG formulations for European, American, Asian, and Barrier options, and their sensitivity. Besides, computational experiments is performed to inspect the numerical efficiency of the method for each option and corresponding Greek. An HPC code for option pricing with the DPG method is provided to motivate the market makers and researchers to utilize the DPG method by customizing the code for their applications and more complicated problems.   
 \bibliography{main}

%% else use the following coding to input the bibitems directly in the
%% TeX file.

% \begin{thebibliography}{00}

% %% \bibitem{label}
% %% Text of bibliographic item

% \bibitem{}

% \end{thebibliography}
\end{document}